\documentclass[10pt]{article}
\input{liemacs10.sty}
\usepackage{hyperref}

\renewcommand{\mlabel}{\label}

\begin{document} 

\title{K\"ahler Geometry, Momentum Maps and Convex Sets}
\author{Karl-Hermann Neeb} 

\maketitle 

\begin{abstract}
{These notes grew out of an expos\'e on M.\ Gromov's paper
``Convex sets and K\"ahler manifolds'' (``Advances in 
Differential Geometry and Topology,'' World Scientific, 1990) 
at the 
DMV-Seminar on ``Combinatorical Convex Geometry and Toric Varieties'' in 
Blaubeuren in April `93. Gromov's paper deals with a proof 
of Alexandrov--Fenchel type 
inequalities and the Brunn--Minkowski inequality for finite dimensional 
compact convex sets and their variants for compact K\"ahler manifolds. 
The emphasis of these notes lies on basic details 
and the techniques from various mathematical areas involved in Gromov's arguments.}
\end{abstract}

\section*{Introduction} \mlabel{sec:0}

These notes grew out of an expos\'e of the author on M.\ Gromov's paper \cite{Gr90} at the 
DMV-Seminar on ``Combinatorical Convex Geometry and Toric Varieties'' in 
Blaubeuren in April `93. Gromov's paper deals with connections 
of {\it Alexandrov--Fenchel 
type inequalities} and the {\it Brunn--Minkowski inequality} (for subsets of 
$\R^n$) with inequalities for the volume of irreducible complex projective 
varieties obtained by B.~Tessier \cite{Te82} and A.~Hovanskii \cite{Ho84}. 

The {\it Brunn--Minkowski inequality} for two compact convex subsets 
$Y_1, Y_2 \subeq \R^n$ asserts that their Minkowski sum 
$Y_1 + Y_2 := \{  a + b \: a \in Y_1, b \in Y_2\}$ 
satisfies 
  \begin{equation}
    \label{eq:4.1x}
\vol(Y_1 + Y_2)^{1\over n} \geq 
\vol(Y_1)^{1\over n} + \vol(Y_2)^{1\over n}.  
  \end{equation}
For bounded convex subsets $Y_1, \ldots, Y_k \subeq \R^n$ and a multiindex 
$J = (j_1,\ldots, j_k) \in \N_0^k$, the {\it mixed volumes} 
$[Y^J] := [Y_1^{j_1}, \ldots, Y_k^{j_k}]$ are defined by the coefficients 
of the homogeneous polynomial 
\[ \vol(t_1 Y_1 + \ldots + t_k Y_k) 
= \sum_{|J| = n} b_J t^J  [Y_1^{j_1}, \ldots, Y_k^{j_k}], \qquad 
b_J := \frac{n!}{j_1! \cdots j_k!},\] 
on $(\R_+)^n$. The {\it Alexandrov--Fenchel Theorem for convex sets} 
asserts that, for $J \in \Delta^n_k$ (the set of all multiindices of degree 
$|J| \leq n$ in $\N_0^k$) the corresponding mixed volumes $[Y^J]$ define a 
function $\Delta_k^n \to \R, J \mapsto \log [Y^J]$ 
which is concave on every ``discrete line'' parallel to the edges. 

In \cite{Gr90} these inequalities are transferred to compact K\"ahler manifolds, 
where one obtains the following 
{\it Brunn--Minkowski inequality for K\"ahler manifolds} 
(Theorem~\ref{thm:IV.12}): 
Let $W_1, \ldots, W_k$ be compact $n$-dimensional 
K\"ahler manifolds and $M \subeq W_1 \times 
\ldots W_k$ be a compact connected 
complex sub\-ma\-ni\-fold of complex dimension $n$ and 
$M_j \subeq W_j$ the projection of $M$ to $W_j$. 
Then 
\begin{equation}
  \label{eq:4.3x}
 \vol(M)^{1\over n} \geq \sum_{j=1}^k \vol(M_j)^{1\over n}. 
\end{equation}
Here one has equality for $n = 1$. 
Replacing the K\"ahler manifolds $W_j$ by projective spaces and $M$ by an irreducible  
projective variety, the inequality remains valid. 
For $n = 2$ this is the Hodge index theorem \cite{GH78}, 
and for $n \geq 3$ and $k = 2$, this inequality is due to Hovanskii and Tessier 
\cite{Te82, Ho84}. 

The central point of \cite{Gr90} is to establish connections between the following 
three contexts: 
\begin{itemize}
\item[\rm(A)] Convex sets, 
\item[\rm(B)] Convex functions, and   
\item[\rm(C)] K\"ahler structures 
\end{itemize}
These connections are so tight that they permit us to translate 
certain theorems such as the Brunn--Minkowski inequality 
which originally belongs to (A) to the areas (B) and (C). The passage from (A) 
to (B) is based on Legendre and Laplace tranforms and representing convex sets 
as the image of the differential of a convex function 
(Fenchel's Convexity Theorem~\ref{thm:I.7}). 
The passage from (B) to (C) is based on K\"ahler potentials on tube domains 
$T_\cD = \cD + i \R^n$ and the quotients $T_\cD/i\Z^n$, where $\cD \subeq \R^n$ 
is an open convex subset, and the corresponding Hamiltonian actions of 
$\R^n$, respectively the torus $\T^n = \R^n/\Z^n$, obtained from translations in the imaginary 
direction. Finally, the passage from (C) back to (A) is established by 
momentum maps for Hamiltonian $\T^n$-actions on toric varieties (which correspond to 
the case $\cD = \R^n$) and their convexity properties as developed in 
\cite{At82, At83, GS82} for compact manifolds and for non-compact manifolds 
and proper momentum maps in \cite{HNP94}. For a generalization from abelian 
to real reductive groups and corresponding gradient maps, see \cite{HS10}. 

Tn the first three sections we describe these three translation mechanisms. 
In Section~\ref{sec:4} we turn to Gromov's results and explain 
how the Brunn--Minkowski inequality is transferred to (B) and (C). 
Eventually we explain 
how it fits into a broader context also including a version of the  
Alexandrov--Fenchel inequality in (C). 
The emphasis of this note lies rather 
on the techniques from various mathematical areas 
involved than on the final results contained in \cite{Gr90}. 
This note does not contain new results, but we hope that it supplies 
useful details in a reasonably self-contained fashion 
and shed some additional light on how 
convex geometry interacts with K\"ahler geometry.\begin{footnote}
{After we wrote the first version of these 
notes in 1993, they were a starting point for 
Chapter V on convex sets and functions in the monograph \cite{Ne00}. So we added 
references to this book in Section~\ref{sec:2}. We also added a last section on 
related subjects and more recent developments in this area. The work 
in \cite{Ne98, Ne99} was also very much inspired by \cite{Gr90}. It constitutes an 
extension of 
the connections between (A), (B) and (C) to the context of non-compact 
non-abelian Lie groups.}
\end{footnote}

\tableofcontents

\section{Convex sets and convex functions} 
\mlabel{sec:1}

In this section we establish a translation mechanism between convex sets and 
convex functions. 

\begin{defn}
Let $V$ be a real finite dimensional vector space.  
A function $f \: V \to \R_\infty := \R \cup \{ +\infty\}$ is called 
{\it convex} if its {\it epigraph}
$$ \epi(f) := \{ (v,t) \in V \times \R \: f(v) \leq t \} $$
is a convex set and the set 
$D_f = f^{-1}(\R)$, the {\it domain} of $f$, is non-empty.
We say that $f$ is {\it closed} if $\epi(f)$ is closed. 

If $C$ is a convex subset of $V$ and $f$ is a convex function on $C$, then we 
think of $f$ as a function $f \: V \to  \R_\infty$ by 
extending $f$ to take the value $\infty$ on the complement of~$C$. 
Note that this does not enlarge the epigraph. 
\end{defn}

\begin{rem} If the convex function $f$ is not closed, then we can always 
consider the closure $\oline{\epi(f)}$ of its epigraph, which is a closed 
convex subset of $V \times \R$. One can show that 
this set is the epigraph of the convex function 
\[ \oline f(x) := \inf \{ t \: (x,t) \in \oline{\epi(f)} \} \] 
(\cite[Prop.~V.3.7]{Ne00}). 
\end{rem}

Let $U \subeq V$ be an open convex set. Further we consider a function 
$f \in C^2(U)$ which has the property that for each $x \in U$ the symmetric 
bilinear form $\dd^2 f(x)$ is positive definite. This implies in 
particular that $f$ is a strictly convex function in $U$ because for each pair 
of different 
points $x,y \in U$ the function $t \mapsto f\big(t x + (1-t)y\big)$ on 
$]0,1[$ has a positive second derivative.

\begin{lem}
\mlabel{lem:I.1}   
The differential $\dd f \: U \to V^*, x \mapsto \dd f(x)$ 
maps $U$ diffeomorphically onto the open subset $\dd f(U)$ of $V^*$. 
\end{lem}

\begin{prf}
It is clear that $\dd f$ is a $C^1$-map. Since the second 
differential $\dd^2 f(x)$ is regular for each $x \in U$, it follows 
from the inverse function theorem that $\dd f$ has a local inverse everywhere. 
We conclude in particular that $\dd f(U)$ is an open set. 

To see that $\dd f$ maps $U$ diffeomorphically onto its image, it remains to prove 
injectivity. This follows from 
$$ \la \dd f(y) - \dd f(x), y - x \ra = \int_0^1 \dd^2 f\big(x + t(y -x)\big)(y-x,y-x)\ 
dt > 0 $$
for $x \not=y$.
\end{prf}

We are interested in properties of open subset $\dd f(U)$ of $V^*$ and in particular in 
those cases where this set is convex. This is not always the case 
(cf.\ \cite{Ef78}), 
so we will have to impose a stronger condition on the function $f$. 
To see what this condition could be, we first have to recall some concepts 
from Fenchel's duality theory for convex functions. 

\begin{defn}
Let $f$ be a convex function on $V$. Then 
\[ f^* \: V^* \to \R_\infty, \qquad f^*(\alpha) := \sup(\alpha - f)  \] 
is called the {\it conjugate} of $f$. 
\end{defn}

The following proposition shows that the passage from $f$ to $f^*$ is similar to the
passage to the adjoint of an unbounded linear operator on a Hilbert space 
(cf.\ \cite{Ru73}). 

\begin{prop}
  \mlabel{prop:I.2} 
{\rm(\cite[Th.\ 12.2]{Ro70})} 
The conjugate $f^*$ of the convex function $f$ on $V$ has the 
following properties:
\begin{itemize}
\item[\rm(i)] $f^*$ is a closed convex function. 
\item[\rm(ii)] $f^*(\alpha) + f(x) \geq \alpha(x)$ for all $x \in V, 
\alpha \in V^*$.
\item[\rm(iii)] $(f^*)^* = \oline f$.
\end{itemize}
\end{prop}

\begin{prf} (i) A pair $(\alpha, t)$ is contained in $\epi(f^*)$ if and only if 
$f^*(\alpha) = \sup(\alpha -f) \leq t$, which is equivalent to 
\begin{equation}
  \label{eq:aff}
(\forall x \in D_f) \quad t - \alpha(x) \geq - f(x).
\end{equation}
Therefore $\epi(f^*)$ is an intersection of a family of closed half spaces 
and therefore a closed convex set. 

(ii) is immediate from the definition of $f^*$. 

(iii) In view of \eqref{eq:aff}, $(\alpha,t) \in \epi(f^*)$ if and only if 
$f \geq \alpha - t$ on $V$, and this is equivalent to 
$\oline f \geq \alpha -t$, so that $f^* = (\oline f)^*$. 

Next we observe that $\oline f \geq \alpha - f^*(\alpha)$ for all $\alpha \in V^*$ 
leads to 
$\oline f \geq f^{**}$. If $\alpha + t \leq \oline f$, then 
$f^*(\alpha) \leq -t$, so that 
$\alpha(x) \leq f^*(\alpha) + f^{**}(x) \leq f^{**}(x) - t$, 
i.e., $\alpha + t \leq f^{**}$. Now the assertion follows from 
\[ \oline f 
= \sup \{ \alpha + t \: \alpha + t \leq f\} 
= \sup \{ \alpha + t \: \alpha + t \leq \oline f\} 
= \sup \{ \alpha + t \: \alpha + t \leq f^{**}\} = f^{**} \] 
(\cite[Lemma~V.3.9(ii)]{Ne00}). 
\end{prf}

\begin{defn}
Let $f \: V \supeq D_f \to \R$ be a convex function on $V$. 
Then every linear functional $\alpha \in V^*$ with 
\[  f(y) \geq f(x) + \alpha(y-x) \qquad \forall y \in D_f,\] 
is called a {\it subgradient of $f$ at $x$}. Note that the preceding 
condition is equivalent to 
\[ f^*(\alpha) = \max(\alpha - f) = \alpha(x) - f(x).\] 
The set of all subgradients of 
$f$ at $x$ is called the {\it subdifferential of $f$ 
at $x$}. It is denoted $\partial f(x)$. It is easy to see that 
$\partial f(x) =\{ \dd f(x)\}$ if $f$ is differentiable in $x$ 
(cf.\ \cite[Th.~25.1]{Ro70}). Note that 
$\partial f(x) = \eset$ if $x \not\in D_f$.
\end{defn}

The following proposition makes the symmetry between $f$ and $f^*$ apparent. 

\begin{prop} \mlabel{prop:I.3} 
{\rm(\cite[Th.\ 23.5]{Ro70}, \cite[Prop.~V.3.22]{Ne00})} Let $f$ be a closed convex function on $V$. Then the 
following are equivalent: 
\begin{itemize}
\item[\rm(1)] $\alpha \in \partial f(x)$.
\item[\rm(2)] $f(x) + f^*(\alpha) = \alpha(x)$.
\item[\rm(3)] $x \in \partial f^*(\alpha)$.
\end{itemize}
\end{prop}

\begin{defn}
(a) In the following we write $\algint(C)$ for the relative interior (also called 
the {\it algebraic interior}) of a convex set 
$C$ with respect to the affine subspace it generates. 

(b) A subset $C$ of a vector space $V$ is said to be {\it almost convex} if 
$\oline C$ is convex and $\algint \oline C$ is contained in $C$. 
\end{defn}

The following observation basically asserts the convexity of the 
range of the multi-valued mapping $\partial f$. It is a first approach to the 
result on the convexity of the image of $\dd f$ for differentiable convex 
functions. 

\begin{cor}
  \mlabel{cor:I.4} If $f$ is a closed convex function, then 
$\im(\partial f)$ is an almost convex dense subset of $D_{f^*}$.
\end{cor}

\begin{prf} (\cite[Cor.~V.3.23]{Ne00}) In view of Proposition~\ref{prop:I.3}, 
$\alpha \in \im(\partial f)$ if and only if 
$\partial f^*(\alpha) \not=\eset$. 
On the other hand, $\partial f^*(\alpha) = \eset$ for $\alpha \not\in D_{f^*}$ and 
$\partial f^*(\alpha) \not= \eset$ for $\alpha \in \algint D_{f^*}$
 (\cite[Th.\ 23.4]{Ro70} or \cite[Lemma~V.3.21(iii)]{Ne00}). 
Therefore $\im(\partial f)$ 
is an almost convex set contained in $D_{f^*}$ which contains $\algint(D_{f^*})$. 
\end{prf}

This corollary shows us that in the case where $f$ is differentiable on the 
interior of its domain, the difficulties are caused by the sets $\partial f(x)$ 
for $x \in \partial U$. So we need a condition which guarantees that these sets 
do not contribute to the image of $\partial f$. 

\begin{defn}
  \mlabel{def:I.5} Let $U \subeq V$ be an open convex subset. 
A function $f\in C^2(U)$ with $\dd^2f(x)$ positive definite for 
all $x \in U$ is called a {\it $C^2$-Legendre function} if for every $x \in U$ and 
every $y \in \partial U$ we have that 
\[  \lim_{t \to 1-} \dd f\big(x + t(y-x)\big)(y-x) = \infty. \] 
\end{defn}

We will see that the $C^2$-Legendre functions are those which are of interest for our 
purposes. 

\begin{lem} {\rm(\cite[Th.\ 26.1]{Ro70} or \cite[Lemma~V.3.30]{Ne00})} \mlabel{lem:I.6}
If $f$ is a closed convex function such that $\Int(D_f) \not= \eset$ 
and the restriction of $f$ to $\Int(D_f)$ is a $C^2$-Legendre function, then 
$\partial f(x)=\eset$ for $x \in \partial D_f$. 
\end{lem}

The following theorem is a sharper version of a result of Fenchel (cf.\ 
\cite{Fe49}). Fenchel did assume that $f$ tends to infinity at the boundary of $U$ 
which is not necessary.

\begin{thm} \mlabel{thm:I.7} 
{\rm(Fenchel's Convexity Theorem
  \begin{footnote}{See \cite{Gr90} for a $C^1$-version of this theorem.}\end{footnote})} Let $U$ 
be an open convex set and $f \: U \to \R$ a $C^2$-Legendre 
function. Then $\dd f(U) = \Int(D_{f^*})$ is an open convex set, 
$\dd f$ maps $U$ diffeomorphically onto $\dd f(U)$, 
$f^*\res_{\Int(D_{f^*})}$ is $C^2$, and $\dd(f^*) \: \dd f(U) \to U$ 
is an inverse of $\dd f$. 
\end{thm}

\begin{prf} We know from Lemma~\ref{lem:I.1} that $\dd f(U)$ is an open set. 
Further, $\oline f$ is a closed convex function, $\epi(\oline f) = \oline{\epi(f)}$, 
$U$ is dense in $D_{\oline f}$, and $\oline f \res_U = f$ (\cite[Prop.~V.3.2(iii)]{Ne00}). 

Now we use Lemma~\ref{lem:I.6}, to see that $\im (\partial\oline f) = \dd f(U)$, so that 
$\dd f(U) = \Int(D_{f^*})$ by Corollary~\ref{cor:I.4} since an open dense almost convex 
subset of $D_{f^*}$ must be equal to $\Int(D_{f^*})$.  

In view of Proposition~\ref{prop:I.3}, it is clear that $\partial f^*(\alpha) 
= \eset$ for $\alpha \not\in \Int(D_{f^*})$, that $\partial f^*$ is single 
valued on $\Int D_{f^*}$, and that it is an inverse of the function 
$\dd f \: U \to \Int(D_{f^*})$. 
\end{prf}

\begin{cor}
  \mlabel{cor:I.8}  Let $f \in C^2(V)$ be such that $\dd^2f(x)$ is positive definite 
for all $x \in V$. Then $\dd f(V) = \Int D_{f^*}$ is an open convex set, 
$\dd f$ maps $V$ diffeomorphically onto $\dd f(V)$, and $\dd(f^*) \: \dd f(V) \to V$ 
is an inverse of $\dd f$. 
\end{cor}

\begin{prf}
This is an immediate consequence of Theorem~\ref{thm:I.7} because in the case 
where $U = V$ the additional condition for a Legendre function is trivially 
satisfied.
\end{prf}

So far we have considered general convex functions. Now we turn to the special 
class of Laplace transforms of (positive) measures. 

\subsection*{Laplace transforms}

Let $V$ be a finite dimensional real vector space and $\mu$ 
a non-zero positive Borel measure on $V^*$. We define the 
{\it Laplace transform} of $\mu$ to be the function 
$$ \cL(\mu) \: V \to ]0,\infty], \qquad x \mapsto \int_{V^*} 
e^{\alpha(x)}\ d\mu(\alpha). $$

\begin{defn}
We define the functions $e_x(\alpha) := e^{\alpha(x)}$ on $V^*$.
Note that $\cL(\mu)(x) > 0$ for all $x \in V$ since 
$\mu \not= 0$. We say that $\mu$ is {\it admissible} if there exists an 
$x \in V$ 
such that $\cL(\mu)(x) < \infty$. If $\mu$ is admissible, then there exists an $x 
\in V$ such that $e_x \mu$ is a finite measure and since $e_x$ is bounded from 
below on every compact set, it follows in particular that $\mu$ is a Radon 
measure and therefore $\sigma$-finite. We write $C_\mu$ for the closed convex 
hull of the support of $\mu$. This is the smallest closed convex subset of 
$V^*$ such that its complement is a $\mu$-null set. 
\end{defn}

\begin{prop}
  \mlabel{prop:I.9}
\begin{itemize}
\item[\rm(i)] The functions $\cL(\mu)$ and $\log \cL(\mu)$ are closed, convex and 
if $C_\mu$ has interior points, then $\cL(\mu)$ is strictly convex on $D_{\cL(\mu)}$.
\item[\rm(ii)] The function $\cL(\mu)$ is analytic on $\Int D_{\cL(\mu)}$ and it has a 
holomorphic extension to the tube domain $\Int D_{\cL(\mu)} + i V$.
\item[\rm(iii)] If $C_\mu$ has interior points, then the bilinear form 
$\dd^2 (\log \cL(\mu))(x)$ is positive definite for all $x \in \Int D_{\cL(\mu)}$. 
\item[\rm(iv)] $D_{\cL(\mu)^*}$ is a convex set which is dense in $C_\mu$.
\end{itemize}
\end{prop}

\begin{prf}
(i) \cite[Th.\ 7.1]{BN78} or \cite[Prop.~V.4.3, Cor.~V.4.4]{Ne00}
\par\nin (ii) \cite[Th.\ 7.2]{BN78} or \cite[Cor.~V.4.4, Prop.~V.4.6]{Ne00}
\par\nin (iii) (\cite[Prop.~V.4.6]{Ne00}) 
In view of (ii), the function $\log \cL(\mu)$ 
is analytic on $\Int D_\mu$. 
We calculate 
$$ \dd(\log \cL(\mu))(x) = {1 \over \cL(\mu)(x)} \dd \cL(\mu)(x) 
= {1 \over \cL(\mu)(x)} \int_{C_\mu} \alpha e^{ \alpha( x) } \ d\mu(\alpha) . $$
Hence 
\begin{align*}
&\ \ \ \ \dd^2(\log \cL(\mu))(x)(y,y) \\
&= -{1 \over \cL(\mu)(x)^2} 
\left(\int_{C_\mu} \alpha(y) e^{\alpha(x) } \ d\mu(\alpha)\right)^2 
 + {1 \over \cL(\mu)(x)} \int_{C_\mu} \alpha(y)^2 e^{\alpha(x)} 
\ d\mu(\alpha) .
\end{align*}
To see that this expression is positive for $0 \not= y$, let 
$g(\alpha) := e^{{1\over 2}\alpha(x)}$
and $h(\alpha) = \alpha(y)g(\alpha)$. 
Then 
\begin{align*}
&\ \ \ \cL(\mu)(x)^2 \dd^2(\log \cL(\mu))(x)(y,y) \\
&= - \left(\int_{C_\mu} g(\alpha) h(\alpha)\ d\mu(\alpha)\right)^2 + 
\int_{C_\mu} g(\alpha)^2 \ d\mu(\alpha) 
\int_{C_\mu} h(\alpha)^2 \ d\mu(\alpha) \geq 0 
\end{align*}
by the Cauchy--Schwartz inequality. In the case of equality it follows that 
$$ \alpha \mapsto h(\alpha) g(\alpha)^{-1} = \alpha(y)$$
is a $\mu$-almost everywhere constant function, i.e., $\supp(\mu)$ lies in an 
affine hyperplane, contradicting the hypothesis that $C_\mu$ 
has interior points. 

\par\nin (iv) \cite[Th.\ 9.1]{BN78} 
\end{prf}

Note that Proposition~\ref{prop:I.9}(iv) shows in particular that the image of 
$\dd \cL(\mu)$
is contained in the convex set $C_\mu$.
The following results gives a good criterion for the convexity of the set
$\dd \cL(\mu)(\Int D_{\cL(\mu)})$: 

\begin{thm}
  \mlabel{thm:I.10} {\rm(Convexity Theorem for Laplace transforms)} 
Suppose that $\Int D_{\cL(\mu)}\not= \eset$ and $\Int C_\mu \not= 
\eset$. Then 
$$ \dd \cL(\mu)(\Int D_{\cL(\mu)}) = \Int C_\mu $$ 
if and only if the restriction of $\cL(\mu)$ to $\Int D_\mu$ is a Legendre 
function. 

If only $\Int D_{\cL(\mu)} \not= \eset$ and $\alpha_0 + H$ is the affine subspace 
generated by $C_\mu$, then $\cL(\mu)$ factors to a function on $V/H^\bot$
and the first assertion applies to the translated measure 
$(-\alpha_0)^*\mu$ on $H \cong (V/H^\bot)^*$. 
\end{thm}

\begin{prf} (\cite[Thm.~V.4.9]{Ne00}) 
The first part follows from \cite[Th.\ 9.2]{BN78}. For the second part we 
note that for $y \in H^\bot$ we clearly have that 
$$ \cL(\mu)(x + y) = \int_{C_\mu} e^{\alpha(x + y)} \ d\mu(\alpha) 
= e^{\alpha_0(y)} \int_{C_\mu} e^{\alpha(x)} \ d\mu(\alpha) 
= e^{\alpha_0(y)} \cL(\mu)(\alpha). $$
Translating $\mu$ by $-\alpha_0$, we obtain a measure $\mu'$ which is supported 
by the subspace $H$, $\cL(\mu')(x) = e^{-\alpha_0(x)} \cL(\mu)(x)$ is its Laplace 
transform which in fact factors to a function on $V/H^\bot$, where it is the 
Laplace transform of $\mu'$, and 
$\dd(\log \cL(\mu'))(x) = -\alpha_0 + \dd(\log \cL(\mu))(x).$ This proves the second 
part since now the first part applies to the measure $\mu'$.
\end{prf}

\begin{cor}
  \mlabel{cor:I.11} If $D_{\cL(\mu)} = V$, 
then $\dd \cL(\mu)(V) = \algint(C_{\mu})$. 
\end{cor}

\begin{thm}
  \mlabel{thm:I.12} 
For every open convex subset $C \subeq V^*$ there exists 
a $C^2$-Legendre function $f$ on $V$ such that $\dd f(V) = C$. 
\end{thm}

\begin{prf} (\cite[Prop.~V.4.14]{Ne00}) 
If $C$ is bounded, then we define $\mu$ to be Lebesgue measure on $C$. 
In general we identify $V$ and $V^*$ with $\R^n$, and define 
$\mu$ by $d\mu(y)= e^{-\|y\|^2} dy$ on $C$ and zero outside of $C$. Then 
$D_{\cL(\mu)} = V$ and $\Int C_\mu = C \not= \eset$. Therefore 
$\dd \cL(\mu)(V)= C$ by Corollary~\ref{cor:I.11}.
\end{prf}

For later applications in Section~\ref{sec:4} we record the following two lemmas. 

\begin{lem}
  \mlabel{lem:I.13} If $f_1$ and $f_2$ are $C^2$-Legendre functions on $V$, then 
\[ \dd(f_1 + f_2)(V) = \dd f_1(V) + \dd f_2(V). \] 
\end{lem}

\begin{prf}
The inclusions 
$$ \dd(f_1 + f_2)(V) \subeq \dd f_1(V) + \dd f_2(V) \subeq 
D_{f_1^*} + D_{f_2^*} \subeq  D_{(f_1 + f_2)^*} $$
follow directly from the definitions and Proposition~\ref{prop:I.3}.  
Since $f_1 + f_2$ is also a 
$C^2$-Legendre 
function on $V$, we conclude that $\dd(f_1 + f_2)(V) = \Int D_{(f_1 + f_2)^*}$ so 
that the left most set is dense in the right most set. 
Now the fact that $\dd f_1(V) + \dd f_2(V)$ is an open convex set proves 
the equality with $\dd(f_1 + f_2)(V)$.
\end{prf}

\begin{lem}
  \mlabel{lem:I.14}Let $U \subeq \R^n$ be an open convex set and $f \in C^2(U)$ such 
that $\dd^2f(x)$ is positive definite for all $x \in U$. Then 
$$ \vol\big(\dd f(U)\big) = \int_{U} \det(\dd^2f)(x) dx, $$
where we identify $\dd^2f$ with the matrix whose entries are 
$\left({\partial^2 f \over \partial x_i x_j}\right)_{i,j = 1,\ldots, n}$. 
\end{lem}

\begin{prf}
Since $\dd^2f(x)$ is the Jacobian of the mapping $\dd f \: U \to \R^n$, and 
it has positive determinant, the assertion follows from the transformation 
theorem for integrals: 
\[  \vol\big(\dd f(U)\big) = \int_{\dd f(U)} 1 \ dx = \int_{U} |\det(\dd^2f(x))|\ dx 
= \int_{U} \det(\dd^2f(x))\ dx. \qedhere\]
\end{prf}

\section{Convex functions and K\"ahler structures}
\mlabel{sec:2}

Now we turn to the relations between convex functions and K\"ahler 
structures. A good reference for foundational material concerning this 
section is Chapter 0.2 in \cite{GH78}. 

We consider $\C^n \cong \R^{2n}$ 
as a $2n$-dimensional real vector space and 
write $J$ for the linear mapping representing the multiplication by $i$. 

\begin{defn}
  \mlabel{def:II.1} (a) A skew symmetric real bilinear form $\omega$ on 
$V = \R^{2n}$ is said to be 
\begin{itemize}
\item[\rm(1)] {\it positive} if $\omega(v,Jv)\geq 0$ holds for all $ v \in V$. 
\item[\rm(2)] {\it strictly positive} if $\omega(v,Jv)> 0$ holds for all 
$ v \in V \setminus \{0\}$. 
\item[\rm(3)] {\it a $(1,1)$-form} if $\omega(Jv,Jw) = \omega(v,w)$ holds for 
all $v,w \in V$.
\end{itemize}

\nin Note that $\omega$ is a $(1,1)$-form if and only if 
$h(v,w) := \omega(v, Jw)$ defines a real symmetric bilinear form on $V$. 

(b) Using the concepts of (a) in each point, we get similar concepts 
for differential $2$-forms on open subsets of $\C^n$ and more generally on 
a complex manifold (one only needs an almost complex structure). 
\end{defn}

Let $M$ be a complex manifold and $J$ the corresponding almost complex 
structure, i.e., for every $x \in M$ the mapping $J_x \: T_x(M) \to T_x(M)$ 
represents multiplication by $i$. Then $J$ acts simply by multiplication on 
vectors, hence on tensors, 
and via duality its also acts on differential forms such that the pairing 
between forms and vectors is invariant under $J$. It follows 
in particular that, if $\omega$ is a $1$-form and ${\cal X}$ a vector field, 
$$ \la J\omega, {\cal X} \ra = \la \omega, J^{-1}{\cal X} \ra = - \la \omega, 
J{\cal X} \ra. $$

Let ${\cal E}^{(p,q)}(M)$ denote the space of $(p,q)$-forms on $M$. Recall 
that a $(p,q)$-form can be expressed in 
local coordinates $z_1,\ldots, z_n$ as a sum
$$ \omega = \sum_{I,J} f_{I,J} \dd z_I \wedge \dd\oline z_J, $$
where 
$I = (i_1, \ldots, i_p) \in \{1,\ldots, n\}^p$ and 
$J = (j_1, \ldots, j_q) \in \{1,\ldots, n\}^q$ 
are multiindices, $\dd z_I = \dd z_{i_1} \wedge \ldots \wedge \dd z_{i_p}$, 
$\dd \oline z_J = \dd \oline z_{j_1} \wedge \ldots \wedge \dd \oline z_{j_q}$, and 
the $f_{I,J}$ are smooth functions. 

Since 
$$J {\partial \over \partial x_j} = {\partial \over \partial y_j} 
\quad \hbox{ and 
} \quad J {\partial \over \partial y_j} = -{\partial \over \partial x_j},$$ 
it follows that $J \dd x_j = \dd y_j$ and $J \dd y_j = - \dd x_j$ which in turn 
yields 
$$ J \dd z_j = -i \dd z_j \quad\hbox{ and } \quad J \dd \oline z_j = i \dd\oline z_j. 
$$
Thus $J \omega = i^{q-p} \omega$ holds for every $(p,q)$-form 
$\omega$ and the $J$-invariant forms are precisely those with $p \equiv  q \mod 4$.

Next we decompose the exterior derivative $\dd \: {\cal E}^m \to 
{\cal E}^{m+1}$ as $\dd = \partial + \oline \partial$, where 
$\partial {\cal E}^{(p,q)} \subeq {\cal E}^{(p+1,q)} $ and 
$\oline\partial {\cal E}^{(p,q)} \subeq {\cal E}^{(p,q+1)}$. On functions we 
have locally 
\begin{equation}
  \label{eq:pq}
\partial f = \sum_j {\partial f \over \partial z_j} \dd z_j, \quad 
\oline\partial f = \sum_j {\partial f \over \partial \oline z_j} \dd\oline z_j
\end{equation}
and on a $(p,q)$-form $\omega = \sum_{I,J} f_{I,J} \dd z_I \wedge 
\dd \oline z_J$ we find 
$$ \partial \omega = \sum_{I,J} \partial f_{I,J} \wedge \dd z_I \wedge \dd \oline 
z_J \quad \hbox{ and } \quad 
\oline\partial \omega = \sum_{I,J} \oline\partial f_{I,J} 
\wedge \dd z_I \wedge \dd \oline z_J. $$

Note that $\dd^2 = 0$ yields $\partial^2 = \oline\partial^2 = \partial \oline \partial 
+ \oline \partial \partial = 0$. 

\begin{lem}
  \mlabel{lem:II.2} Let $f$ be a smooth function on the complex manifold $M$. Then 
the following assertions hold: 
\begin{itemize}
\item[\rm(i)] $ \dd J \dd f = 2 i \partial \oline \partial f $
is a $(1,1)$-form and in particular $J$-invariant. 
\item[\rm(ii)] If $h_f$ is the symmetric form defined by 
$$ h_f(v,w) := \dd J\dd f(v,Jw), $$
then 
$h_f = \dd^2 f + J \dd^2 f$ in local coordinates. 
\end{itemize}
\end{lem}

\begin{prf}
(i) We write $\dd f = \partial f + \oline \partial f$. Then 
$$J \dd f = J \partial f + J \oline \partial f = -i \partial f + i \oline\partial f
= i(\oline \partial - \partial)f. $$
Therefore 
$$ \dd J \dd f = i(\partial + \oline \partial)(\oline\partial-\partial)f = 
2i \partial \oline\partial f. $$
Since $\oline\partial f$ is a $(0,1)$-form, $\partial \oline\partial f$ is a 
$(1,1)$-form and the assertion follows. 

\par\nin (ii) First we calculate 
$$ \dd J\dd f = 2i \sum_{j,k} {\partial^2 f \over \partial z_j \partial \oline z_k} 
\dd z_j \wedge \dd\oline z_k. $$
Writing 
$$ \dd z_j \wedge \dd\oline z_k = \dd z_j \otimes \dd\oline z_k - \dd\oline z_k \otimes 
\dd z_j, $$
we find with \eqref{eq:pq} that 
\begin{align*}
h_f 
&= 2i \sum_{j,k} {\partial^2 f \over \partial z_j \partial \oline z_k} 
(\dd z_j \otimes J^\top \dd\oline z_k - \dd\oline z_k \otimes J^\top \dd z_j) \\
&= 2i \sum_{j,k} {\partial^2 f \over \partial z_j \partial \oline z_k} 
(-i \dd z_j \otimes \dd\oline z_k - i \dd\oline z_k \otimes \dd z_j) \\
&= 2\sum_{j,k} {\partial^2 f \over \partial z_j \partial \oline z_k} 
(\dd z_j \otimes \dd\oline z_k + \dd\oline z_k \otimes \dd z_j). 
\end{align*}
On the other hand we have 
\begin{align*}
 \dd^2 f 
&= \sum_{j,k} \Big({\partial^2 f \over \partial z_j \partial z_k} 
\dd z_j \otimes \dd z_k + {\partial^2 f \over \partial z_j \partial \oline z_k} 
(\dd z_j \otimes \dd\oline z_k + \dd\oline z_k \otimes \dd z_j) 
 + {\partial^2 f \over \partial \oline z_j \partial \oline z_k} 
\dd \oline z_j \otimes \dd\oline z_k\Big). 
\end{align*}
In $\dd^2 f + J \dd^2 f$, the first and the last term cancel and we obtain 
$\dd^2 f + J \dd^2 f = h_f$. 
\end{prf}

Note that the definition of $H$ rests on the choice of local coordinates but 
that the definition of $h_f$ is coordinate free. 

We pave the way to later applications with the following proposition 
which relates Legendre functions to exact positive $(1,1)$-forms on complex 
tori.

\begin{prop}
  \mlabel{prop:II.3}
Let $M = \C^n / i \Z^n \cong (\C^\times)^n$ be endowed with the usual complex structure.  
\begin{itemize}
\item[\rm(i)] Let $U \subeq \R^n$ be an open subset and $f \in C^\infty(U)$ 
be such that $\dd^2f(x)$ is positive 
semidefinite for all $x \in U$. We define the function $\tilde f$ on 
$M_U := \{ z + i \Z^n \: \Re z \in U\}$ by 
$\tilde f(z + i\Z^n) := f(\Re z)$. Then 
$\omega_f := \dd J \dd\tilde f$ is a positive 
$(1,1)$-form which is strictly positive if and 
only if $\dd^2 f$ is everywhere positive definite. Moreover 
$$ {1\over n!} \int_{M_U} \omega_f^n = \int_U \det(\dd^2 f). $$
\item[\rm(ii)] If $\omega$ is a positive exact $(1,1)$-form on $M$ which is 
invariant under the action of the torus $T = i\R^n/i\Z^n$, then there 
exists a smooth function $f$ on $\R^n$ such that $\omega = \omega_f$ and 
$\dd^2f$ is everywhere positive semidefinite. 
\end{itemize}
\end{prop}

\begin{prf}
(i) We have already seen in Lemma~\ref{lem:II.2} that $\omega :=\omega_f$ is a 
$(1,1)$-form. 
Let $v,w \in \R^n$ and $m = x + i y + i\Z^n \in M$. Then we use 
Lemma~\ref{lem:II.2}(ii) 
to see that 
\begin{align*}
&\ \ \  \omega(m)\big(v + i w, i(v + iw)\big)
= \omega(m)(v + i w, -w + iv) \\
&= \dd^2 \tilde f(m)(v + i w, v + i w) + \dd^2 \tilde f(m)(w - i v, w- iv) \\
&= \dd^2 f(x)(v, v) + \dd^2 f(x)(w, w) \geq 0. 
\end{align*}
This calculation also shows that $\omega$ is strictly positive if and only if 
$\dd^2f$ is everywhere positive definite. 

To prove the integral formula, we first use the above calculation and 
polarization to obtain 
$$ \omega(m)(v + i w, v' + iw') = \dd^2f(x)(v,w') - \dd^2f(x)(w,v'). $$
Thus we can write 
\begin{equation}
  \label{eq:2.1}
\omega = \sum_{j,k} { \partial^2 f \over \partial x_j \partial x_k} \dd x_j 
\wedge \dd y_k.  
\end{equation}
Let $\alpha_k := \sum_{j=1}^n { \partial^2 f \over \partial x_j \partial x_k} \dd x_j$. 
Then $\omega = \sum_k \alpha_k \wedge \dd y_k$ and therefore 
\begin{align*}
 {1\over n!} \omega^n 
&= \alpha_1 \wedge \dd y_1 \wedge \alpha_2 \wedge \ldots \wedge \dd y_n 
= \sum_{j_1, \ldots, j_n=1}^n 
{\partial^2 f \over \partial x_{j_1} \partial x_1} \dd x_{j_1} \wedge \dd y_1 \wedge 
{\partial^2 f \over \partial x_{j_2} \partial x_2} \dd x_{j_2} \wedge \ldots
\wedge \dd y_n \\
&= \det(\dd^2 f) \dd x_1 \wedge \dd y_1 \wedge \dd x_2 \wedge \ldots \wedge \dd y_n. 
\end{align*}
Now we calculate the integral: 
\begin{align*}
{1\over n!} \int_M \omega^n 
&= \int_{\big(\R \times (\R/\Z)\big)^n} 
\det(\dd^2 f)(x) d x_1 d y_1 d x_2 d y_2 \ldots d x_n d y_n
= \int_{\R^n} \det(\dd^2 f)(x) d x. 
\end{align*}

(ii) First we observe that every invariant exact $2$-form 
$\beta = \dd\alpha$ on the torus $T$ is also the differential of an invariant exact form, 
hence zero. This can be achieved by averaging $\alpha$. On the other hand every invariant 
$1$-form on $T$ is closed because $T$ is abelian (cf.\ \cite{ChE48}). 
Therefore $T$ permits no non-zero invariant exact $2$-forms.  

This shows that the restriction of $\omega$ to the $T$-cosets in 
$M$ vanishes. Since $\omega$ is also $J$-invariant, it follows that 
$i\R^n$ is also everywhere isotropic for $\omega$. Therefore, in coordinates 
from the mapping $\R^{2n} \to M$, we have that 
$$ \omega(x) = \sum_{j,k} a_{jk}(x) \dd x_j \wedge \dd y_k, $$
where we have used the $T$-invariance to see that the functions 
$a_{jk}$ do not depend on the $y_k$'s. 

The closedness of $\omega$:
$$ \dd\omega = \sum_{j,k,\ell} {\partial a_{jk} \over \partial x_\ell} \dd x_\ell \wedge 
\dd x_j \wedge \dd y_k = 0 $$
now shows that the $1$-forms $\alpha_k := \sum_{j,k} a_{jk} \dd x_j$ are also 
closed. Hence there exist smooth functions $f_k$ on $\R^n$ with 
$\dd f_k = \alpha_k$, i.e., $a_{jk} = {\partial f_k \over \partial x_j}$. 

The $J$-invariance of $\omega$ further yields 
$$ \omega = J \omega = \sum_{j,k} a_{jk} \dd y_j \wedge (-\dd x_k) 
= \sum_{j,k} a_{jk} \dd x_k  \wedge \dd y_j $$
and therefore $a_{jk} = a_{kj}$. Thus  
$$ {\partial f_k \over \partial x_j} = {\partial f_j \over \partial x_k} $$
which means that $\sum_k f_k \dd x_k$ is closed. Therefore we find a smooth 
function $f$ on $\R^n$ with $\dd f = \sum_k f_k \dd x_k$. Finally $\omega = \dd J \dd 
\tilde f$ follows from \eqref{eq:2.1}. 
\end{prf}

\section{K\"ahler structures, Hamiltonian actions 
and momentum maps}
\mlabel{sec:3}

In this section we describe how the momentum map for Hamiltonian group actions 
can be used to relate K\"ahler manifolds and convex sets. We apply this in 
particular to orbits of complex tori in projective spaces. Again we do not 
intend to describe the geodesic way to the final results but rather to 
show how the different approaches in the literature fit together. 

For more details on symplectic manifolds, momentum maps and Hamiltonian 
actions we refer to \cite{GS84} and \cite{LM87}.
\begin{defn}
  \mlabel{def:III.1}  (a) Let $M$ be a smooth real manifold. A {\it symplectic 
structure} on $M$ is a closed, non-degenerate $2$-form $\omega$. 

(b) Let $M$ be a complex manifold. A {\it K\"ahler structure} on 
$M$ is a strictly positive, closed $(1,1)$-form $\omega$. 

Then $\omega$ is in particular non-degenerate, so that it 
defines a symplectic structure on $M$. If $J$ denotes the almost 
complex structure of $M$, then for each $x \in M$ the sesquilinear form 
$$ h(x)(v,w) := \omega(x)(v, Jw) - i \omega(x)(v,w) $$
is positive definite, so that it defines a complex Hilbert space structure on 
the tangent space $T_x(M)$. 
\end{defn}

Let $(M,\omega)$ be a connected 
symplectic manifold. Then we can associate to each
smooth function $F \in C^\infty(M)$ a {\it Hamiltonian vector field} 
${\cal X}_F$ which is determined by the equation 
$$ \dd F = - i_{{\cal X}_F} \omega. $$
If we define a Lie algebra structure on $C^\infty(M)$ by the {\it Poisson 
bracket} 
$$ \{ F, G \} := \omega({\cal X}_F, {\cal X}_G) =  {\cal X}_F G = - {\cal X}_G 
F, $$
then 
$$ C^\infty(M) \to {\cal V}(M), \quad F \mapsto {\cal X}_F $$
is a homomorphism of Lie algebras whose kernel consists of the constant 
functions. 

Now let $G$ be a connected Lie group acting on $(M,\omega)$ and leaving the 
symplectic structure invariant. Then we obtain a homomorphism of Lie algebras 
$\psi \: \g = \Lie(G) \to {\cal V}(M)$ with 
\[ \psi(X)_p :=  \frac{d}{dt}\Big|_{t=0} \exp(-tX).p \quad \mbox{ for all } \quad 
p \in M.\]
 Such an action is called {\it Hamiltonian} if there
exists a homomorphism $\phi \: \g \to C^\infty(M)$ in such a way that 
$$ {\cal X}_{\phi(X)} = \psi(X) $$
holds for all $X \in M$, i.e., $\psi$ lifts to a homomorphism 
$\g \to (C^\infty(M), \{\cdot,\cdot\})$. 

For a Hamiltonian action one has a mapping 
$$ \Phi \: M  \to \g^*, \quad \Phi(m)(X) := \phi(X)(m)\big)$$
called the {\it momentum map}. It is an equivariant map 
from $M$ to $\g^*$, where $G$ acts on $\g^*$ by the coadjoint action 
$g.\alpha := \Ad^*(g)\alpha := \alpha \circ \Ad(g)^{-1}$.

In the last section we have seen how to obtain K\"ahler structures on a complex 
manifold $M$ 
via smooth functions as $\omega = \dd J \dd f$. Now we bring this together 
with holomorphic Hamiltonian actions of a connected Lie group $G$ on $M$.

Let $G \times M \to M$ be an action of $G$ on $M$ by holomorphic mappings 
and $\psi \: \g \to {\cal V}(M)$ the corresponding homomorphism into the Lie 
algebra of vector fields. For $X \in \g$ we set 
\begin{equation}
  \label{eq:dagg}
\phi(X) := \la J \dd f, \psi(X) \ra = - \big(J \psi(X)\big)f \quad \hbox{ for } 
f \in C^\infty(M). 
\end{equation}

\begin{prop}
  \mlabel{prop:III.2}Suppose that the function $\psi(X)f$ is constant for every 
$X \in \g$. Then $\Phi(m)(X) := \phi(X)(m)$ defines a momentum map for the 
action of $G$ on the symplectic manifold $(M,\omega)$.
\end{prop}

\begin{prf}
Let $X \in \g$ and write ${\cal L}_{\psi(X)}$ for the Lie derivative 
along the vector field $\psi(X)$. 
Then the fact that $G$ acts holomorphically means that ${\cal L}_{\psi(X)} J = 
0$ for all $X \in \g$ and therefore 
$$ {\cal L}_{\psi(X)} J \dd f 
= J {\cal L}_{\psi(X)}\dd f 
= J \dd i_{\psi(X)} \dd f = J \dd \big(\psi(X)f\big) = 0 $$
by the Cartan formula. Hence 
$$ 0 = {\cal L}_{\psi(X)} J \dd f
= i_{\psi(X)} \dd J \dd f + \dd i_{\psi(X)}J \dd f 
= i_{\psi(X)} \omega + \dd \phi(X) $$
shows that ${\cal X}_{\phi(X)} =\psi(X)$. 

It remains to show that $\phi$ is a homomorphism of Lie algebras. To this end, 
we first note that 
$J [\psi(X), \psi(Y)] = [J \psi(X), \psi(Y)]$ follows from the fact that $G$ 
acts holomorphically because the Lie algebra of all vector fields 
generating holomorphic flows is a complex Lie subalgebra with respect to the complex 
structure induced by~$J$. 
We also recall that $\psi(Y)f$ is constant so that 
${\cal X}\psi(Y)f = 0$ holds for every vector field ${\cal X}$ on $M$.
With these remarks in mind we calculate: 
\begin{align*}
\phi([X,Y]) 
&= - \big(J\psi([X,Y])\big)f 
= - \big(J[\psi(X),\psi(Y)]\big)f 
= - \big([J\psi(X),\psi(Y)]\big)f \\
&= \psi(Y)\big(J \psi(X)\big)f - J\psi(X)\big(\psi(Y)f\big) 
= \psi(Y)\big(J \psi(X)\big)f 
= -\psi(Y) \phi(X) \\
&= -{\cal X}_{\phi(Y)} \phi(X) 
= -\{ \phi(Y), \phi(X)\} = \{ \phi(X), \phi(Y)\}. 
\qedhere\end{align*}
\end{prf}

\begin{rem}
  \mlabel{rem:III.3} If 
the vector field $\psi(X)$ on $M$ has compact orbits, then 
the constancy of the function $\psi(X)f$ implies that even $\psi(X)f = 0$ 
because every periodic affine function is constant. It follows in particular 
that, for a compact group $G$, the function $f$ has to be invariant under the action of $G$.
\end{rem}

\subsection*{Projective spaces}

A particular interesting case for these constructions is the complex projective 
space $\bP(\C^{n})$. For a non-zero vector $v \in \C^{n}$ we write 
$[v] = \C^\times v$ for the corresponding ray in $\bP(\C^{n})$ and 
write $(v_1 : \ldots : v_n)$ for $[(v_1, \ldots, v_n)]$ 
(homogeneous coordinates). 

We want to define a K\"ahler structure $\omega$ on the projective space
(the {\it Fubini--Study metric}) (cf.\ \cite[p.~30]{GH78}). 
The most transparent way to do this is 
to construct $\pi^*\omega$ on the unit sphere $\bS^{2n-1}$, where 
$\pi \: \C^n\setminus \{0\} \to \bP(\C^n)$ is the canonical projection. 

We start with the function $F(z) := {1\over 2}\log \|z\|^2$. Then 
$$ \dd F(z)(v) = {1\over \|z\|^2} \Re \la z,v \ra, \quad 
J \dd F(z)(v) = {1\over \|z\|^2} \Im \la v,z \ra, $$
and 
$$ \dd J \dd F(z)(v,w) 
= { 2 \over \|z\|^4}\left( \Re \la z,v\ra\Im \la z,w \ra 
- \Re \la z,w\ra\Im \la z,v \ra\right)  
+ { 2 \over \|z\|^2}\Im \la w,v \ra. $$
For $w = J v$, this specializes to 
$$ \dd J \dd  F(z)(v,Jv) 
= -{ 2 |\la z,v\ra|^2 \over \|z\|^4}   
+ { 2 \|v\|^2 \over \|z\|^2}
= { 2 \over \|z\|^2}\left(\|v\|^2 - \left|\left\la {z\over 
\|z\|},v\right\ra\right|^2\right). $$
This number only depends on the length of the projection of $v$ to the complex 
hyperplane orthogonal to the ray $\C z$. 

Therefore there exists a smooth $2$-form $\omega$ on $\bP(\C^n)$ such that 
$\pi^*\omega\res_{\bS^{2n-1}} = \dd J\dd F$. 
It follows in particular that $\omega$ is $J$-invariant, 
hence a $(1,1)$-form, and that it is strictly positive. Since the unitary group 
$\U_n(\C)$ acts by holomorphic maps on $\C^n$ and $F$ is constant under this action, 
the action of $\U_n(\C)$ preserves the K\"ahler 
structure on $\bP(\C^n)$. 

Often it is more convenient to have the form $\omega$ in homogeneous coordinates. 
So let $z \in \C^n$ with $\|z\| = 1$. We obtain a chart 
$\beta \: z^\bot \to \bP(\C^n)$, $v \mapsto [v + z]$. We want to calculate the 
pullback $\beta^*\omega$ as a $(1,1)$-form on the $(n-1)$-dimensional complex 
vector space $z^\bot$. Since we know already that the pullback to the unit 
sphere is $\dd J \dd F$, we simply have to calculate the pullback of 
$\dd J\dd F$ under the mapping $\beta' \: v \mapsto {1\over \|v + z\|}(v + z)$. 
Let $x , v,w \in z^\bot$. Then 
$$ \dd\beta'(x)(v) = {1\over \|x + z\|} \left( v - { x + z \over \|x + z\|^2}  
\Re \la x + z,v \ra\right). $$
Therefore, using that the form $\dd J\dd F$ depends in $\beta'(x)$ 
only on the or\-tho\-go\-nal projection of the vectors on the complex hyperplane 
$\beta'(x)^\bot$ and that 
$\dd J\dd F(\lambda z) = {1\over \|\lambda\|^2} \dd J\dd F(z)$, we see that  
\begin{align*}
\dd J\dd F\big(\beta'(x)\big)(\dd\beta'(x)v,\dd\beta'(x)w)
&= {1\over \|x + z\|^2} \dd J\dd F\left({ x + z \over \|x+ z\|}\right)(v,w)\\
&= \dd J\dd F(x + z)(v,w). 
\end{align*}
Hence the pullback $\beta^*\omega$ is simply the restriction of $\dd J\dd F$ 
to the hyperplane $z + z^\bot$ in $\C^n$, but this is the form $\dd J\dd F_z$, where 
$$F_z(x) = F(x + z) = {1\over 2}\log \|x+z\|^2 
= {1\over 2}\log(1 + \|x\|^2).$$ 

\msk
Here we are merely interested in linear torus actions on projective space, so 
we consider the following situation. Let $V \cong \C^k$ be a finite dimensional 
complex Hilbert 
space, $\bP(V)$ the projective space endowed  with the Fubini--Study metric, 
$ T = \R^n/\Z^n$, $T_{\C} = (\C^\times)^n \cong \C^n/i\Z^n$, and $\pi \: 
T_{\C} \to \GL(V)$ a holomorphic representation. 
Then we find an orthonormal basis $e_1, \ldots, e_k$ in $V$ (which we use to 
identify $V$ with $\C^k$), and holomorphic characters 
$\chi_1, \ldots, \chi_k$ of $(\C^\times)^n$ such that 
$\pi(z)e_j = \chi_j(z) e_j$ for $j = 1,\ldots, k$ and $z \in T_{\C}$. 
We write $\chi_j(z) = e^{-i\alpha_j(z)}$ for $\alpha_j \in \ft^*$, where $\ft = 
\Lie(T)$. 

Let $\Psi \: \C^{k-1} \to \bP(\C^k)$ be a coordinate 
chart obtained by homogeneous 
coordinates by $\Psi(v') = (1 : v_2 : \ldots : v_k)$, 
where $v' = (v_2, \ldots, v_k)$. 
Then the image of $\Psi$ in the projective space is invariant under the induced 
action of the group $T_{\C}$. In homogenous coordinates the action of 
$T_{\C}$
is given by 
$$ z.(1 : v_2 : \ldots : v_k) = \big(1 : \chi_1(z)^{-1}\chi_2(z)v_2 :  \ldots 
: \chi_1(z)^{-1}\chi_k(z) v_k \big). $$
Let $\pi' \: T_{\C} \to \GL_{k-1}(\C)$ denote the representation defined by the 
characters $\chi_2\chi_1^{-1}, \ldots, \chi_k \chi_1^{-1}$. 
Then the function $F_{e_1}(v') = {1\over 2}\log(1 + \|v'\|^2)$ is 
invariant under the action of $T$. So we can use Proposition~\ref{prop:III.2} 
to see that 
the action of $T$ is Hamiltonian and with $v = (1,v_2, \ldots, v_k)$ the momentum 
map is given by 
\begin{eqnarray} \label{eq:3.1} 
\Phi\big([(1,v')]\big)(X) 
&=& -\big(J\psi(X)F)_{e_1}(v') = \dd F_{e_1}(v')\big(-i\psi(X)v'\big) 
= {1 \over 1 + \|v'\|^2} \Re \la v', i Xv'\ra \notag\\
&=& i {\la Xv', v'\ra \over 1 + \|v'\|^2} 
= \sum_{j=2}^k (\alpha_j-\alpha_1)(X){|v_j|^2 \over 1 + \|v'\|^2} \notag\\
&=& -\alpha_1(X) { \|v'\|^2 \over  1 + \|v'\|^2} + 
\sum_{j=2}^k \alpha_j(X){|v_j|^2 \over 1 + \|v'\|^2} \notag\\
&=& -\alpha_1(X) + \sum_{j=1}^k \alpha_j(X){|v_j|^2 \over \|v\|^2}. 
\end{eqnarray}
Therefore 
$$ \Phi\big([(1,v')]\big) 
= \sum_{j=2}^k (\alpha_j -\alpha_1) { |v_j|^2 \over 1 + \|v'\|^2}. $$

The next step is to investigate the image of a $T_{\C}$-orbit in projective 
space under the momentum map. Let $z = x + i y + i \Z^n \in T_{\C}$. Then 
\begin{align*}
\Phi\big(z.[(1,v')]\big) 
&= { \sum_{j=2}^k (\alpha_j - \alpha_1) |e^{-i(\alpha_j - \alpha_1)(z)} v_j|^2 
\over  1 + \sum_{j=2}^k |e^{-i(\alpha_j - \alpha_1)(z)} v_j|^2 } 
= { \sum_{j=2}^n (\alpha_j - \alpha_1) e^{2(\alpha_j - \alpha_1)(y)} |v_j|^2 
\over  1 + \sum_{j=2}^k e^{2(\alpha_j - \alpha_1)(y)} |v_j|^2 }.
\end{align*}

Write $\delta_\alpha$ for the Dirac measure concentrated in the point $\alpha 
\in \ft^*$ and let 
\[ \mu = \delta_0 + \sum_{j=2}^k |v_j|^2 \delta_{\alpha_j - \alpha_1}
= \sum_{j=1}^k |v_j|^2 \delta_{\alpha_j - \alpha_1}.\] 
Then the above formula 
for the momentum mapping shows that 
\begin{equation}
  \label{eq:ddagg}
\Phi\big(z.[(1,v')]\big) = \dd(\log \cL(\mu))(2y), 
\end{equation}
where $\cL(\mu)$ is the Laplace transform of the measure $\mu$ 
(cf.\ Section~\ref{sec:1}). 
Now the results of Section~\ref{sec:1} make it easy for us to compute the image of the 
$T_{\C}$-orbit of the line $[(1,v')]$ as the relative 
interior of the convex hull of 
the support of $\mu$ (Corollary~\ref{cor:I.11}). Therefore 
\begin{equation}
  \label{eq:3.2}
\Phi\big(T_{\C}.[(1,v')]\big) 
= \algint(\conv\{ \alpha_j \: j = 1,\ldots, k\}) -\alpha_1. 
\end{equation}

This result shows that the natural momentum map which was defined by the chart 
obtained by homogeneous coordinates depends on the choice of this chart. 
Of course, since $T$ is abelian, we could easily define another momentum map by 
taking $\phi'(X) := \phi(X) + \gamma(X)\1$, where 
$\1$ is the function constant 
to $1$ and $\gamma$ is a linear functional on $\ft$. Then the image of the momentum 
map is shifted by $\gamma$. In view of formula \eqref{eq:3.2}, it seems to be natural 
to take $\gamma = \alpha_1$. That this choice is in fact a rather natural one 
can be seen as follows. 

We have already observed that the Fubini--Study metric on $\bP(\C^k)$ is 
invariant under the action of the unitary group $\UU_k(\C)$. We will see that the
action of $\U_k(\C)$ is in fact Hamiltonian. Writing $\fu_k(\C)$ for the Lie algebra 
of $\U_k(\C)$ which consists of the space of skew-Hermitean matrices, 
we define for each $X \in \fu_k(\C)$ the function 
$$ \phi'(X)([v]) := i { \la Xv,v \ra \over \la v,v\ra} $$
on $\bP(\C^k)$ (note that the right hand side only depends on $[v] = \C^\times v$). 
As we have seen above in homogeneous coordinates, 
this is up to a constant the correct Hamiltonian function, hence 
$i_{{\cal X}_{\phi'(X)} }\omega + \dd\phi'(X) = 0$. On the other hand 
\begin{align*}
\{\phi'(X), \phi'(Y)\}([v]) 
&= {\cal X}_{\phi'(X)}\phi'(Y)([v]) 
= \frac{d}{dt}\Big|_{t=0} \phi'(Y)([e^{-tX}v]) 
= \frac{d}{dt}\Big|_{t=0} i { \la Y e^{-tX}v, e^{-tX}v \ra \over \la e^{-tX}v, 
e^{-tX}v \ra} \\
&= \frac{d}{dt}\Big|_{t=0} i { \la e^{tX} Y e^{-tX}v, v \ra \over \la v, v \ra} 
= i { \la [X,Y]v, v \ra \over \la v, v \ra} 
= \phi'([X,Y])([v]). 
\end{align*}
Therefore $\Phi'([v])(X) := \phi'(X)([v])$ defines a momentum map for the action 
of $\U_k(\C)$ on $\bP(\C^k)$. 

If $\pi \: T_{\C} \to \GL_k((\C)$ is a representation of the complex torus 
given by the characters $\chi_1, \ldots , \chi_k$ as above, we find in 
particular that 
$$ \Phi'([v])(X) = i { \sum_{j=1}^k (-i\alpha_j)(X) |v_j|^2
\over \la v,v\ra} $$
and therefore 
$$ \Phi'([v]) = {1\over \|v\|^2} \sum_{j=1}^k |v_j|^2 \alpha_j. $$
This shows that the image of the whole projective space is the 
convex hull of the set $P :=\{ \alpha_1, \ldots, \alpha_k\}$ of all weights. 
But as we have already seen, 
much stronger results hold. To make this 
explicit with the new momentum map, let $z = x + i y + i \Z^n \in T_{\C}$. 
Then 
\begin{equation}
  \label{eq:3.3}
\Phi'(z.[v]) 
= { \sum_{j=1}^k \alpha_j e^{2\alpha_j(y)}|v_j|^2 \over 
\sum_{j=1}^k e^{2\alpha_j(y)}|v_j|^2} = \dd(\log f_{\mu_v})(2y) 
\quad \mbox{ for } \quad \mu_v = \sum_{j=1}^k |v_j|^2 \delta_{\alpha_j}. 
\end{equation}

\begin{thm}
  \mlabel{thm:III.4} 
Let $T_{\C} = \C^n/i\Z^n$ and $\pi$ be a holomorphic 
representation on $\C^n$ given by the characters $\chi_1, \ldots, \chi_k$ 
with $\chi_j(z) = e^{-i\alpha_j(z)}$ and $\alpha_j \in \ft^*$. 
Let further $[v] \in \bP(\C^k)$ and consider the holomorphic 
action of $T_{\C}$ on $\bP(\C^k)$ induced by $\pi$. Then the momentum map 
$$ \Phi'([v]) = {1\over \|v\|^2} \sum_{j=1}^k |v_j|^2 \alpha_j $$
for the Hamiltonian action of the torus $T = i\R^n/i\Z^n$ on $\bP(\C^k)$ 
satisfies
\begin{equation}
  \label{eq:3.4}
\Phi'(T_{\C}.[v]) 
= \algint(P_v), \quad \mbox{ where } \quad 
P_v := \conv\{\alpha_j \: v_j \not=0\}
\end{equation}
is a polyhedron in $\ft^*$. 

The orbit closure $\oline{T_{\C}.[v]}$ is mapped onto the whole polyhedron 
$P_v$. The other $T_{\C}$-orbits are mapped onto the relative interior of the 
faces of $P_v$. This establishes a bijection between $T_{\C}$-orbits in the 
orbit closure and faces of $P_v$. 
\end{thm}

\begin{prf}
The first assertion about the image of $T_{\C}.[v]$ follows from \eqref{eq:3.3} 
and Corollary~\ref{cor:I.11}. 
Since the orbit closure is compact, it has a compact image 
and since the orbit is dense (the Zariski closure and the closure in the 
manifold topology coincide) and $\Phi'$ is continuous, it follows that 
$\Phi'$ maps it onto $P_v$. 

Now let $[v']$ be contained in the closure of $T_{\C}.[v]$. Then, by taking 
another element in $T_{\C}.[v']$ if necessary, we can assume that there exists 
a set $J \subeq \{1,\ldots, k\}$ such that 
\[  v_j' =
\begin{cases}
v_j & \text{for } j \in J  \\
0  & \text{ for } j \not\in J  
\end{cases}\]
and $F_J := \conv \{ \alpha_j : j \in J\}$ is a face of $P_v$ (cf.\ 
\cite[Th.\ IV.13]{Ne92}, \cite[Prop.\ 1.6]{Od88}, or \cite[p.416]{Od91}). 
Then it follows from the first part of the theorem that 
$$ \Phi'(T_{\C}.[v']) = \algint(F_J)$$
and that two different orbits are mapped into different faces. 
\end{prf}

Orbit closures as we have seen above are called {\it projective toric 
varieties}. Therefore Theorem~\ref{thm:III.4} establishes a connection between 
the projective toric variety $\oline{T_{\C}.[v]}$ and the 
polyhedron in $\ft^*$ which arises as the image of the momentum map and which 
describes the stratification of the orbit closure into orbits. 

More information on toric varieties and momentum mappings can be found in 
\cite{Ju81, At82, At83, GS82, Br85}. For a discussion of the normality 
of the orbits closures $\oline{T_{\C}.[v]}$ we refer to \cite[pp.95/96]{Od88}.

\section{Mixed volumes and inequalities}
\mlabel{sec:4}

In this section we eventually turn to the applications of the 
techniques explained above to the Brunn--Minkowski inequality 
and the Alexandrov--Fenchel inequality. 

In the following we fix $k \in \N$ and 
write $I \subeq \N_0^k$ for a multiindex $I = (i_1, \ldots, i_k)$. 
We define its {\it degree} by $|I| := \sum_{j=1}^k i_j \in \N_0$ and  write 
$\Delta^n_k$ for the set of all $k$-multiindices of degree $n$ and 
$\Delta^{\leq n}_k$ for the set of all $k$-multiindices with degree 
$\leq n$. We consider $\Delta_k^{\leq n}$ as a ``discrete simplex" of dimension~$k$.

\begin{lem}
  \mlabel{lem:IV.1}Let $n \in \N$, $\K$ a field of characteristic $0$, 
and for a multiindex $I$ write 
$Iy := \sum_{j=1}^k i_j y_j$ for the corresponding linear form on $\K^k$. 
Then the polynomials $(Iy)^n$, $I \in \Delta^n_{k}$ 
form a basis of the space of homogeneous 
polynomials of degree $n$ in $y_1, \ldots, y_k$. 
\end{lem}

\begin{prf}
For a multiindex $J \in \Delta_{k}^n$ we write 
$$ b_J := { n! \over j_1! \ldots j_k!} $$
for the corresponding binomial coeffficient. Then 
$$ (Iy)^n = \sum_{J \in \Delta_{k}^n} b_J (i_1 y_1)^{j_1} \cdot \ldots \cdot (i_k 
y_k)^{j_k} = \sum_{J \in \Delta_{k}^n} b_J I^J y^J. $$
Therefore the square 
matrix with entries $a_{IJ} = b_J I^J$ describes the coefficients of the 
polynomials $(Iy)^n$ in the monomial basis. To show that this matrix is 
regular, it suffices to deal with the matrix given by $a'_{IJ} = I^J$. 
We thus have to show that $\sum_J c_J I^J =0$ for all $I \in \Delta_{k}^n$ 
implies $c_J = 0$ for all $J$. 

This is equivalent to show that, for a homogeneous polynomial 
$f(y) = \sum_{|J|=n} c_J y^J$ to vanish, it suffices to vanish on 
$\Delta_{k}^n$. We consider the affine map 
$$ \phi \: \K^{k-1} \to \K^k, \quad (y_1, \ldots, y_{k-1}) \mapsto 
(y_1, \ldots, y_{k-1}, n - y_1 - \ldots - y_{k-1}) $$
which maps $\Delta_{k-1}^{\leq n}$ onto $\Delta_k^n$. 
Therefore $f$ vanishes on $\Delta_k^n$ if and only if the inhomogeneous polynomial 
$g := f \circ \phi$ of degree $\leq k$ vanishes on $\Delta_{k-1}^{\leq n}$. 
We thus have to show that a polynomial $g(y) = \sum_{J \in \Delta_{k-1}^{\leq n}} 
c_J y^J$ of $k-1$ variables $y_1, \ldots, y_{k-1}$
vanishes if it vanishes on the set $\Delta_{k-1}^{\leq n}$. This means that 
the subset $\Delta_{k-1}^{\leq n} \subeq \K^{k-1}$ is a {\it determining subset} 
for the space of polynomials of degree $\leq k-1$. 

This is done by induction. Since it is trivial for $k = 1$ 
or $n = 1$ (affine maps), 
we assume that $k,n > 1$ and that the assertion is true for all 
smaller values of $k$ and $n$. 
Suppose that $g(I) = 0$ for all $I \in \Delta_{k-1}^{\leq n}$. 
Restricting to the hyperplanes given by $y_\ell = 0$, it follows from our 
induction hypothesis that $c_J = 0$ if $j_\ell = 0$ for some $\ell$. Therefore
$c_J \not= 0$ implies that there exists a multiindex $J'$ with 
$J = J' + (1,1,\ldots, 1)$. Hence 
$$ g(y) = \sum_{J \in \Delta_{k-1}^{\leq n}} c_J y^J = y_1 \cdot\ldots\cdot 
y_{k-1} \sum_{J' \in \Delta_{k-1}^{\leq n-k+1}} 
c_{J' + (1,1,\ldots,1)} y^{J'}. $$
We conclude that the polynomial  $g'(y) = \sum_{J' \in \Delta_{k-1}^{\leq n-k+1}} 
c_{J' + (1,1,\ldots,1)} y^{J'}$
vanishes on the set $(1,\ldots, 1) + \Delta_{k-1}^{\leq n-k+1}$. 
By our induction hypothesis, the set $\Delta_{k-1}^{\leq n-k+1}$ is a determining 
subset for the space of polynomials of degree $\leq n-k+1$. 
Hence all translates of this set are determining, and this leads to 
$g' = 0$, so that also $c_{J' + (1,1,\ldots,1)} = 0$ 
for all $J'\in \Delta_{k-1}^{\leq n-k+1}$, 
and therefore that $c_J = 0$ for all $J \in \Delta_{k-1}^{\leq n}$. 
This completes the 
proof.
\end{prf}

The preceding lemma will permit us later to define mixed volumes via positive 
$(1,1)$-forms. Let $k \in \N$ 
and $\Omega := (\omega_1, \ldots, \omega_k)$ be a sequence of 
positive $(1,1)$-forms on the compact complex manifold $M$ of complex dimension 
$n$. For $I \in \Delta_{k}^n$ we put 
$$ \Omega^I := \omega_1^{i_1} \wedge \ldots \wedge \omega_k^{i_k}. $$
We are interested in the behaviour of the function 
$I \mapsto \int_M \Omega^I$.

\begin{defn}
  \mlabel{def:IV.2} A function $\ell \: \Delta_k^n \to \R$ 
is called {\it $1$-concave} if it is 
concave on every (discrete) line (isomorphic to some $\Delta_1^{\leq m}$) 
parallel to the edges. 

If, for example, $k = 2$, then $\Delta_2^m = \{(0,m), (1,m-1), \ldots, (m,0)\}$
and the {\it $1$-concavity} means that 
\[ \ell\Big(\sum_j \alpha_j (j,m-j)\Big) \geq \sum_j \alpha_j \ell(j,m-j) \] 
whenever $\alpha_j \geq 0$ with $\sum_j \alpha_j = 1$ and 
$\sum_j \alpha_j (j,m-j) \in \Delta_2^m$. 
\end{defn}

\begin{lem}
  \mlabel{lem:IV.3}To check that a function $\ell \: \Delta^n_{k} \to \R$ 
is concave, it suffices to check that 
$$ \ell(y_0) \geq {1\over 2} \ell(y_-) + {1\over 2} \ell(y_+), $$
where $(y_-,y_0,y_+)$ forms a part of a ``discrete line segment" in 
$\Delta^n_{k}$ which is parallel to an edge, i.e., 
$y_+ - y_0 = y_0  - y_- = e_{j_1} - e_{j_2}$ for two canonical 
basis vectors $e_{j_1}$, $e_{j_2}$ of $\R^n$. 
\end{lem}

\begin{prf}
  Since we only have to consider a line segment parallel to an edge of 
$\Delta_{k}^n$, we may w.l.o.g.\ assume that 
$k = 2$. Let $I = \conv \Delta^n_2 \subeq \R^2$. Then we extend $\ell$ to a 
continuous piecewise affine function on $I$. The graph of this function is a
polygon and the condition imposed on $\ell$ yields that this polygon is concave at 
every vertex. Therefore it is a concave polygon and this implies that 
$\ell$ is a concave function.
\end{prf}

It is the preceding lemma which is responsable for the fact that one gets 
merely $1$-concavity in the following theorem.

\begin{thm}
  \mlabel{thm:IV.4} {\rm(Alexandrov--Fenchel Theorem 
for compact complex manifolds)} 
Let $M$ be a connected compact complex manifold, $n = \dim_{\C}M$, and 
$\Omega = (\omega_1, \ldots, \omega_k)$ a sequence of positive closed 
$(1,1)$-forms. Then the function 
$$ \ell : \Delta^n_{k} \to \R^+, \quad I \mapsto \log \int_M \Omega^I $$
is $1$-concave.
\end{thm}

\begin{prf}
For a proof we refer to \cite[Thm.~I.6B]{Gr90}.
\end{prf}

\begin{rem}
  \mlabel{rem:IV.5} 
The preceding theorem remains true for irreducible projective 
varieties in the sense that one has to consider only those 
$(1,1)$-forms $\omega$ on 
the set $M_{\rm reg}$ of regular points which have the property that for every 
holomorphic map $\beta \: U \to M$, $U \subeq \C^n$ open, the pull-back $\beta^*\omega$ on 
$\beta^{-1}(M_{\rm reg})$ extends smoothly 
to $U$. Again we refer to \cite{Gr90}. The main ingredient in the proof is 
Hironaka's theorem on the resolution of singularities (cf.\ \cite{Hi70}).
\end{rem}

\subsection*{Applications to convex sets} 

\begin{defn} \mlabel{def:IV.6} 
Let $Y_1, \ldots, Y_k$ be compact convex subsets of 
$\R^n$. Consider Legendre functions $f_j$ with $\algint(Y_j)= \dd f_j(\R^n)$, 
and write $\omega_j := \omega_{f_j} = \dd J \dd \tilde f_j$ for the corresponding positive 
$(1,1)$-form on the complex manifold $M = \C^n/i\Z^n \cong T_{\C}$ 
(cf.~Proposition~\ref{prop:II.3}).

For $I \in \Delta_{k}^n$ the integral ${1\over n!}\int_M \Omega^I$ 
is called the {\it $I$-th mixed volume} of 
$(Y_1, \ldots, Y_k)$, denoted 
\[  [Y^I]  = [Y_1^{i_1}, \ldots, Y_k^{i_k}] 
= \frac{1}{n!} \int_M \omega_1^{i_1} \wedge \cdots \wedge \omega_k^{i_k}. \] 
\end{defn}

This definition is motivated by the following observation. 
For $t_1, \ldots, t_k \geq 0$ with \break {$\sum_j t_j > 0$}, the function 
$f := t_1 f_1 + \cdots + t_k f_k$ is a Legendre function with 
$\dd f(\R^n) = \algint(Y)$ for $Y = t_1 Y_1 + \cdots + t_k Y_k$ (Lemma~\ref{lem:I.13}). 
By an easy induction we derive from Lemmas~\ref{lem:I.13}, \ref{lem:I.14} and 
Proposition~\ref{prop:II.3}(i) that 
\begin{align*}
 \vol(t_1 Y_1 + \ldots + t_k Y_k) 
&= \vol\big( \dd(t_1 f_1 + \cdots + t_k f_k)(\R^n)\big) \\ 
&= \int_{\R^n} \det\big(d^2(t_1 f_1 + \cdots + t_k f_k)\big)(x)\, dx 
 \qquad \mbox{ by Lemma~\ref{lem:I.14}} \\ 
&= \frac{1}{n!} \int_M \omega_f^n  \qquad \mbox{ by Proposition~\ref{prop:II.3}} \\ 
&= {1\over n!} \int_M (t_1 \omega_1 + \ldots + t_k \omega_k)^n  \\
&= {1\over n!} \int_M \sum_{|J| = n} b_J t^J \omega_1^{j_1} \wedge \cdots \wedge 
\omega_k^{j_k}  \\
&= \sum_{|J| = n} b_J t^J  [Y_1^{j_1}, \ldots, Y_k^{j_k}]
\end{align*}
because $\omega_f = \dd J \dd \tilde f = t_1 \omega_{1} + \cdots + t_k \omega_{k}$. 
This calculation shows that the function 
\begin{equation}
  \label{eq:mixedvol}
\R_+^n \to \R, \quad (t_1, \ldots, t_n) \mapsto 
 \vol(t_1 Y_1 + \ldots + t_k Y_k) = \sum_{|J| = n} b_J t^J  [Y_1^{j_1}, \ldots, Y_k^{j_k}]
\end{equation}
is a polynomial of degree $\leq n$ whose coefficients are determined by the 
mixed volumes. This justifies the terminology. For 
$k = 1$, we obtain in particular 
\begin{align*}
[Y^n] = \vol(Y).
\end{align*}

\begin{lem}
  \mlabel{lem:IV.7}
Let $Y_j$, $f_j$, and $\omega_j$, $j =1,\ldots, k$, be as above. 
Then the following assertions hold: 
\begin{itemize}
\item[\rm(i)] For $I, J \in \Delta_{k}^n$, let $c_{IJ}$ 
be the coefficients for which 
$y^I = \sum_{J} c_{IJ} (Iy)^n.$
Then we have the relation 
$$ \Omega^I = \sum_{J} c_{IJ} (j_1 \omega_1 + \ldots + j_k \omega_k)^n. $$
\item[\rm(ii)] 
$ [Y^I] = \sum_{J} c_{IJ} \vol(j_1 Y_1 + \ldots + j_k Y_k). $
\item[\rm(iii)] $\vol(Y_1 + Y_2) = [(Y_1 + Y_2)^n] = \sum_{j=0}^n {n \choose j} 
[Y_1^j, Y_2^{n-j}].$ 
\end{itemize}
\end{lem}

\begin{prf}
(i) First we note that the existence of the $c_{IJ}$ follows from 
Lemma~\ref{lem:IV.1}. Since the $2$-forms $\omega_j$ generate a commutative algebra, it 
follows from $y^I = \sum_{J} c_{IJ} (Iy)^n$ that 
\[ \Omega^I = \sum_{J} c_{IJ} (i_1 \omega_1 + \ldots + i_k \omega_k)^n.\] 

(ii) For $f := j_1 f_1 + \cdots + j_k f_k$, we have seen above that 
\[  \vol(j_1 Y_1 + \ldots + j_k Y_k) 
= {1\over n!} \int_M (j_1 \omega_1 + \ldots + j_k \omega_k)^n.  \] 
Therefore (ii) follows from (i). 

(iii) is a special case of \eqref{eq:mixedvol}. 
\end{prf}
 
Note that Lemma~\ref{lem:IV.7}(ii) shows in 
particular that $[Y^I]$ does not depend on the choice of the functions $f_j$. 

\begin{thm}
  \mlabel{thm:IV.8}{\rm(Alexandrov--Fenchel Theorem for convex sets)} Let 
$Y_1, \ldots, Y_k$ be bounded convex subsets of $\R^n$. Then 
$$ \ell_Y(I) := \log [Y^I] $$
defines a $1$-concave function on $\Delta^n_{k}$. 
\end{thm}

\begin{prf}
Since one can approximate bounded convex sets arbitrarily well by 
polyhedra generated by rational points, and since $[Y^I]$ is positively 
homogeneous, 
it suffices to prove Theorem~\ref{thm:IV.8} 
for polyhedra with integral extreme points. 
We may even assume that one extreme point of each polyhedron is the origin.

Suppose that $Y_1 = \conv\{ 0,\alpha_1^1, \ldots, \alpha_{\ell_1}^1\}$ with 
$\alpha_j^1(\Z^n) \subeq 2\pi \Z$ for all $j$. 
Then we identify $\R^n$ with $\ft^*$ for $T = i\R^n/i\Z^n$ and consider 
the representation of $T_{\C}$ defined by the characters $\chi_0^1(z) = 1$, 
$\chi^1_j(z) := e^{-i\alpha^1_j(z)}$, $j =1,\ldots, \ell_1$ on $\C^{\ell_1+1}$. 
Let $v_1 = (1,\ldots,1) \in \C^{\ell_1+1}$ and consider the orbit closure 
$M_1 := \oline{T_{\C}.[v_1]}$ in the projective space $\bP(\C^{\ell_1+1})$. 
Then the corresponding momentum map maps $M_1$ onto $Y_1 \subeq \ft^*$ 
(Theorem~\ref{thm:III.4}) in such a way that the pull-back $\tilde \omega_1$ 
of the K\"ahler form $\omega_1$ on $M_1$ via the orbit 
map $T_{\C} \to M_1, z \mapsto z.[v_1]$ can be written as 
$\tilde \omega_1 = \dd J\dd F_1$, where 
\[  F_1(z) = {1\over 2}\log(1 + \|z v_1'\|^2) 
= {1\over 2}\log\Big(1 + \sum_{j=1}^{\ell_1} |\chi_j^1(z)|^2\Big) 
= {1\over 2}\log\Big(1 + \sum_{j=1}^{\ell_1} e^{2\alpha_j^1(y)}\Big).\]
So $F_1 = -\log \cL(\mu)(-2y)$, where $\cL(\mu)$ is the Laplace transform of the 
measure 
$$ \mu = \delta_0 + \sum_{j=1}^{\ell_1} \delta_{\alpha_j}. $$
Then $\dd F_1(y) = \dd(\log \cL(\mu))(2y)$ and therefore $\tilde\omega_1$ 
is a positive $(1,1)$-form on $T_{\C}$ representing the polyhedron $Y_1$. 
We proceed similarly for $Y_2, \ldots, Y_k$. 

Now, by definition of the mixed volume,  
$$ [Y^I] 
= {1\over n!} \int_{T_{\C}} (i_1 \tilde \omega_1 
+ \ldots + i_k \tilde \omega_k)^n . $$
We consider the mapping 
$$ \beta \: T_{\C} \to \bP(\C^{\ell_1+1}) \times 
\ldots \times  \bP(\C^{\ell_k+1}), \quad z \mapsto (z.[v_1], \ldots, z.[v_k]) $$
and the toric variety $M := \oline{\beta(T_{\C})}$. The
projection $\pi_j$ onto the $j$-th factor maps 
$M$ onto $M_j$ and the pull-back $\pi_j^*\omega_j$ satisfies 
$\tilde\omega_j = \beta^*\pi_j^*\omega_j$. Therefore 

\begin{align*}
  [Y^I] 
&= {1\over n!} \int_{T_{\C}} (i_1 \tilde \omega_1 
+ \ldots + i_k \tilde \omega_k)^n  
= {1\over n!} \int_{T_{\C}} \beta^*(i_1 \pi_1^*\omega_1 
+ \ldots + i_k \pi_k^*\omega_k)^n  \\
&= {1\over n!} \int_{M} (i_1 \pi_1^*\omega_1 
+ \ldots + i_k \pi_k^*\omega_k)^n  
= {1\over n!} \int_{M} \Omega^I 
\end{align*}
for $\Omega = (\pi_1^*\omega_1 , \ldots, \pi_k^*\omega_k)$.  
Now Remark~\ref{rem:IV.5} tells us that Theorem~\ref{thm:IV.4} applies and this completes the 
proof. 
\end{prf}

\begin{lem} \mlabel{lem:5.9} 
If $f \: \Delta_k^n \to \R$ is a $1$-concave function, then 
\[ f(i_1, \ldots, i_k) \geq \sum_{j = 1}^k \frac{i_j}{n} f(n e_j) \quad \mbox{ and } \quad 
e^{f(i_1, \ldots, i_k)} \geq \prod_{j = 1}^k (e^{f(ne_j)})^{\frac{i_j}{n}},  \] 
where $e_1, \ldots, e_k \in \Delta_k^1$ are the extreme points. 
\end{lem}

\begin{prf} This is verified by induction on $k$. F\"ur $k \leq 2$ the assertion 
follows from the definition. So let us assume that $k > 2$ and that it holds for 
$k -1$. If one entry $i_j$ vanishes, then the induction hypothesis applies directly. 
If this is not the case, then we observe that, 
for $m := i_1 + i_2$, the element 
$(i_1, \ldots, i_k)$ lies on the discrete line between 
\[ (m, 0, i_3, \ldots, i_k) \quad \mbox{ and } \quad (0,m, i_3, \ldots, i_k).\] 
Therefore $1$-concavity and the induction hypothesis leads to 
\begin{align*}
 f(i_1, \ldots, i_k) 
&\geq \frac{i_1}{m} f(m, 0,i_3, \ldots, i_k) +  
\frac{i_2}{m} f(0, m ,i_3, \ldots, i_k) \\
&\geq \frac{i_1}{m} \Big(\frac{m}{n} f(n e_1) + \sum_{j > 2} \frac{i_j}{n} f(n e_j)\Big)
+ \frac{i_2}{m} \Big(\frac{m}{n} f(n e_2) + \sum_{j > 2} \frac{i_j}{n} f(n e_j)\Big)\\
&= \frac{i_1}{n} f(ne_1) + \frac{i_2}{n} f(ne_2) + 
\frac{i_1 + i_2}{m} \sum_{j > 2} \frac{i_j}{n} f(n e_j)
= \sum_{j = 1}^k \frac{i_j}{n} f(n e_j). \qedhere
\end{align*}
\end{prf}

With Lemma~\ref{lem:5.9} we derive from Theorem~\ref{thm:IV.8}: 

\begin{cor}
  \mlabel{cor:IV.9} $[Y^I] \geq \vol(Y_1)^{i_1 \over n} \cdot \ldots \cdot 
\vol(Y_k)^{i_k \over n}.$
\end{cor}

\begin{prf}
  $\ell_Y(I) = \ell_Y(i_1, \ldots, i_k) \geq \sum_j {i_k \over n} \ell_Y(0,\ldots, 
0,n,0,\ldots,0)$. 
\end{prf}
 
\begin{cor}
  \mlabel{cor:IV.10}  {\rm(Brunn--Minkowski inequality)}
  \begin{equation}
    \label{eq:4.1}
\vol(Y_1 + Y_2)^{1\over n} \geq 
\vol(Y_1)^{1\over n} + \vol(Y_2)^{1\over n}. 
  \end{equation}
\end{cor}

\begin{prf}
We calculate with Lemma~\ref{lem:IV.7}(iii) and Corollary~\ref{cor:IV.9}:
\begin{align*}
\vol(Y_1 + Y_2) 
&= \sum_j {n \choose j} [Y_1^j, Y_2^{n-j}] 
\geq \sum_j {n \choose j} \vol(Y_1)^{j\over n}\vol(Y_2)^{n-j\over n}\\
&= \left( \vol(Y_1)^{1\over n} + \vol(Y_2)^{1\over n}\right)^n.
\qedhere\end{align*}
\end{prf}

\begin{rem}
  \mlabel{rem:IV.11} 
Let $f_1$ and $f_2$ 
be $C^2$-Legendre functions on $\R^n$ and 
$C_{f_j} = \dd f_j(\R^n)$. Then Lemma~\ref{lem:I.14} shows that 
$$ \vol(C_{f_j}) = \int_{\R^n} \det \dd^2f(x) \ dx. $$
Since $C_{f_1 + f_2} = C_{f_1} + C_{f_2}$ by Lemma~\ref{lem:I.13}, we see that 
on the level of convex functions the Brunn--Minkowski inequality reads
\begin{equation}
  \label{eq:4.2}
\left( \int_{\R^n} \det \dd^2(f_1 + f_2)(x)\ dx \right)^{1\over n} 
\geq \left( \int_{\R^n} \det \dd^2 f_1(x)\ dx \right)^{1\over n} 
+ \left( \int_{\R^n} \det \dd^2 f_2(x)\ dx \right)^{1\over n}.
\end{equation}
That, conversely, this inequality also implies \eqref{eq:4.1} 
follows from Theorem~\ref{thm:I.12} which permits us to represent each bounded open 
convex set as $C_f$ for a 
$C^2$-Legendre function $f$ on $\R^n$. 

Using Proposition~\ref{prop:II.3}, we can translate \eqref{eq:4.2} into to the inequality 
$$ \left(\int_{T_{\C}} (\omega_1 + \omega_2)^n \right)^{1\over n} 
\geq \left(\int_{T_{\C}} \omega_1^n \right)^{1\over n} 
\left(\int_{T_{\C}} \omega_2^n \right)^{1\over n} $$
for all positive exact $T$-invariant $(1,1)$-forms on $T_{\C}$. 
\end{rem}

In the context of K\"ahler manifolds, one has the following version of 
the Brunn--Minkowski inequality. 

\begin{thm}
  \mlabel{thm:IV.12} {\rm(Brunn--Minkowski inequality for K\"ahler manifolds)} 
Let $W_1, \ldots, W_k$ be compact connected $n$-dimensional 
K\"ahler manifolds and $M \subeq W := W_1 \times 
\ldots W_k$ be a compact connected complex sub\-ma\-ni\-fold of complex dimension $n$, 
$\pi_j \: M \to W_j$ the projections, 
and $M_j := \pi_j(M)$. Then 
\begin{equation}
  \label{eq:4.3}
 \vol(M)^{1\over n} \geq \sum_{j=1}^k \vol(M_j)^{1\over n}. 
\end{equation}
\end{thm}

\begin{prf}
Let $\omega_j$ denote the K\"ahler form on $W_j$ and $\Omega := 
(\pi_1^*\omega_1, \ldots, \pi_k^*\omega_k)$. Then 
$ \omega := \sum_{j=1}^k \pi^*\omega_j$ defines the induced K\"ahler structure 
on $M$. Now we have 
\begin{align*}
\vol(M)
&= {1\over n!} \int_M \omega^n 
= {1\over n!} \int_M \left(\sum_{j=1}^k \pi^*\omega_j\right)^n \\
&= {1\over n!} \int_M \sum_{|I| = n} b_I \Omega^I  = {1\over n!} \sum_{|I| = n} b_I \int_M \Omega^I \\
&\geq {1\over n!} \sum_{|I| = n} b_I 
\left(\int_M (\pi^*\omega_1)^n\right)^{i_1\over n}  \cdot\ldots \cdot 
\left(\int_M (\pi^*\omega_k)^n\right)^{i_k\over n} 
\qquad \mbox{ by Corollary~\ref{cor:IV.9}} \\
&= {1\over n!} \left(\sum_{j=1}^k \left(\int_M (\pi^*\omega_j)^n\right)^{1\over 
n}\right)^n 
= \left(\sum_{j=1}^k \left({1\over n!} 
\int_{M_j} \omega_j^n\right)^{1\over n}\right)^n 
= \left(\sum_{j=1}^k \vol(M_j)^{1\over n}\right)^n. 
\qedhere\end{align*}
\end{prf}

Note that one has equality for $n = 1$ in the preceding theorem because one 
trivially has equality in Theorem~\ref{thm:IV.4} in this case.

\begin{rem}
  \mlabel{rem:IV.12} 
The preceding result remains true if we replace the 
K\"ahler manifolds $W_j$ by projective spaces and $M$ by an irreducible  
projective variety (cf.\ Remark~\ref{rem:IV.5} and \cite{Gr90}).

For $k = 2$, we have $M \subeq P_1 \times P_2$ and 
$$ \vol(M)^{1\over n} \geq \vol(M_1)^{1\over n}+ \vol(M_2)^{1\over n}. $$

As already mentioned above, we have equality for $n = 1$. 
As explained in \cite[\S 3.3]{Gr90}, for $n =2$ this inequality 
can be derived from Hodge's inequality, resp., Hodge's index theorem in the form 
$$ \left(\int_M \omega_1 \wedge \omega_2\right)^2 \geq 
\left(\int_M \omega_1^2\right) \left(\int_M \omega_2^2\right), $$ 
(cf.\ \cite{GH78}), 
and for $n \geq 3$, this inequality is due to Hovanskii--Tessier. 
Therefore it is called the {\it Hodge--Tessier--Hovanskii inequality} (cf.\ 
\cite{Te82, Ho84}). 
\end{rem}

For further connections of the topics of this paper with algebraic geometry, 
cohomology and K\"ahler manifolds we refer to \cite{Gr90} and 
\cite[pp.102--104]{Od88}. 

\section{Perspectives}
\mlabel{sec:5} 

In this final section we comment on some more recent developments related to the themes 
of \cite{Gr90}.

{\bf Pushforwards of measure by gradients of convex functions:} 
One interesting issue we touched in Section~\ref{sec:1} is writing an open convex 
subset $C \subeq \R^n$ as $\dd f(\R^n)$ for a $C^2$-Legendre function, which 
may be a Laplace transform $f = \cL(\mu)$, where 
$d\mu(y) = e^{-\|y\|^2}dy$ on $C$ (Theorem~\ref{thm:I.12}). A slightly 
different issue is to consider measures on a convex set $C$ which 
are the push-forward $\mu_\psi$ 
under the differential $\dd \psi$ of a convex function 
$\psi \: \R^n \to \R \cup \infty$, of the measure $e^{-\psi(x)}\, dx$. 
Since $\psi$ is locally Lipschitz, its differential $\dd \psi$ exists almost
everywhere, so that such measures make sense for general convex functions. 
Here the finiteness and non-triviality of the measure $e^{-\psi(x)}\, dx$ 
is equivalent to $0 < \int_{\R^n} e^{-\psi(x)}\, dx < \infty$, which 
in turn is equivalent to the domain $D_\psi$ having interior points and 
$\lim_{x \to \infty} \psi(x) = \infty$. In \cite[Thm.~2]{CK15} a class of convex functions 
is determined for which the assignment $\psi \mapsto \mu_\psi$ leads to a 
bijection onto the class of finite Borel measures whose barycenter is the origin 
and whose support spans the whole space. \\

{\bf Optimal transport:}
This connects to optimal transport theory as follows. For a given measure 
$\mu = \mu_\psi$, the differential $\dd \psi$ is the {\it quadratic optimal map}, 
or {\it Brenier map}, between the measure $e^{-\psi(x)}\, dx$ on $\R^n$ and 
the measure $\mu$  (\cite{Br91}). This refers to the existence of the 
unique {\it polar factorization} 
of a measurable map $u \: X \to \R^n$ from a probability space 
$(X,\mu)$ to $\R^n$ in the form $u(x) = \dd\psi(s(x))$, where 
$\Omega \subeq \R^n$ is a bounded domain endowed with the normalized 
Lebesgue measure $\mu_\Omega$, $s \: (X,\mu)  \to (\Omega,\mu_\Omega)$ 
is measure preserving, and $\psi \: \Omega \to \R$ is convex (\cite{Br91}).
Note that, for $\psi(x) = \frac{1}{2}\|x\|^2$, 
the measure $e^{-\psi(x)}\, dx$ simply is the Gaussian measure on $\R^n$. 
For applications of the measures $\mu_\psi$ and momentum maps 
to Poincar\'e type inequalities in analysis, we refer to \cite{Kl13}. \\

{\bf Generalizations and applications of the Brunn--Minkowski inequality (BMI):} 
The BMI for bounded subsets $Y_0, Y_1 \subeq \R^n$ implies
 with $a^2 + b^2 \geq 2ab$ for $a,b \geq 0$ the inequality 
\[ \vol\Big(\frac{Y_0 + Y_1}{2}\Big) \geq \vol(Y_0)^{1\over 2} \vol(Y_1)^{1\over 2}.  \]
This in turn implies that, for $Y_t := (1-t) Y_0 + t Y_1$, the function 
$\log(\vol(Y_t))$ is concave. This observation can be generalized as follows. 
In $\R^{n+1}\cong \R \times \R^n$, we consider the convex subset 
\[ Y := \bigcup_{0 \leq t \leq 1} \{t\} \times Y_t.\] 
Then the function $t \mapsto \vol(Y_t)$ is a marginal of Lebesgue measure 
restricted to $Y$. The $\log$-concavity of this function is a special case 
of {\it Pr\'ekopa's Theorem} asserting that, for any 
convex function $F \: \R^{n+1} \to \R \cup \{\infty\}$ with 
$\int_{\R^{n+1}} e^{-F} < \infty$, the function 
\[ t \mapsto \log \int_{\R^n} e^{-F(t,\bx)}\, d\bx \] 
is concave. In fact, if $F_Y$ is the convex indicator function of $Y$ 
which is $0$ on $Y$ and $\infty$ elsewhere, then the right hand side specializes 
to $\log(\vol(Y_t))$ (see \cite{CK12}). 

For the connections of the BMI and its generalizations in 
various branches of mathematics, we refer to Gardner's nice survey 
\cite{Ga02}. Here one finds in particular a discussion of the BMI for non-convex 
subsets, the situations where equality holds,
 and its applications to isoperimetric inequalities and estimates in analysis. 
Moreover, generalizations of the BMI to the sphere, hyperbolic space, 
Minkowski space and Gauss space (euclidean space where the volume is measured with respect to a Gaussian density) are explained. 
For subsets $A,B$ of the integral lattices $\Z^n$, 
analogs of the BMI giving lower bounds of the 
cardinality of $|A +B|$ in terms of $|A|$ and $|B|$ 
can be found in \cite{GG01}. 

Connections to the representation theory of reductive groups have been 
established by the work of V.~Okounkov who used the BMI to study weight 
polytopes (\cite{Ok96}).
\\

{\bf Infinite dimensional convex geometry:} 
There exist natural infinite dimensional contexts in which substantial portions 
of the duality theory for convex functions work. 
Here one may start with a real bilinear duality pairing 
$(\cdot,\cdot) \: V \times W \to \R$ of two infinite dimensional real vector spaces 
$V$ and~$W$. This pairing defines natural (weak) locally convex topologies on 
$V$ and $W$. Accordingly, one may consider convex functions 
$f \:  V \to \R \cup \{\infty\}$, define closedness in terms of the 
weak closedness of the epigraph $\epi(f) \subeq V \times \R$ and 
consider the conjugate convex function $f^* \: W \to \R \cup\{\infty\}, 
f^*(w) := \sup_{v \in V} (v,w) - f(v)$ which is automatically closed. 
However, for more refined applications, it is important to also have a finer 
locally convex topology on $V$ for which the domain of $f$ has interior points. 
For more details, further developments and applications see 
\cite{Mi08}, \cite{Bou07} and \cite[\S 3]{Ro74}.


\begin{thebibliography}{NRW01}

\bibitem[At82]{At82} Atiyah, M. F., {\it Convexity and commuting hamiltonians}, Bull.
London Math. Soc. {\bf 14} (1982), 1--15

\bibitem[At83]{At83} ---, {\it Angular momentum, convex polyhedra and algebraic geometry}, 
J. of Edinburg Math. Soc. {\bf 26}(1983), 121--138

\bibitem[BN78]{BN78} Barndorff--Nielsen, O., ``Information and Exponential Families -- In 
Statistical Theory", John Wiley and Sons, New York, 1978

\bibitem[Bou07]{Bou07} Bourbaki, N., ``Espaces vectoriels topologiques. 
Chap.1 \`a 5'', Springer-Verlag, Berlin, 2007 

\bibitem[Br91]{Br91} Brenier, Y., {\it Polar factorization and monotone rearrangement 
of vector-valued functions}, Comm. Pure. Appl. Math. {\bf 44:4} (1991), 375--417

\bibitem[Br85]{Br85} Brugi\`eres, A., {\it Propri\'et\'es de convexit\'e de l'application 
moment [d'apr\`es Atiyah, Guillemin Sternberg, Kirwan et. al.]}, Ast\'erisque, 
Seminaire Bourbaki, Exp.\ 654, 1985/86

\bibitem[CK12]{CK12} Cordero-Erausquin, D., and B. Klartag, {\it 
Interpolations, convexity and geometric inequalities}, in 
``Geometric Aspects of Functional Analysis,'' 151--168, 
Lecture Notes in Math. {\bf 2050}, Springer, Heidelberg, 2012

\bibitem[CK15]{CK15} ---, {\it Moment measures}, 
J. Funct. Anal. {\bf 268:12} (2015), 3834--3866

\bibitem[ChE48]{ChE48} Chevalley, C. and S. Eilenberg, {\it Cohomology theory of Lie groups and Lie 
algebras}, Transactions of the Amer. Math. Soc. {\bf 63} (1948), 85--124 
 
\bibitem[Ef78]{Ef78} Efron, B., {\it The geometry of exponential families}, 
Annals of Stat. {\bf 6}(1978), 362--376

\bibitem[Fe49]{Fe49} Fenchel, M., {\it On conjugate convex functions}, Canad. J. Math. {\bf 
1}(1949), 73--77

\bibitem[Ga02]{Ga02} Gardner, R. J., {\it The Brunn--Minkowsi Inequality}, 
Bulletin of the Amer. Math. Soc. {\bf 39:3} (2002), 355--405 

\bibitem[GG01]{GG01} Gardner, R. J., and P.~Gronchi, 
 {\it A Brunn--Minkowsi inequality for the integer lattice}, 
Trans. Amer. Math. Soc. {\bf 353:10} (2001), 3995--4024 

\bibitem[GH78]{GH78} Griffiths, P., and J. Harris, ``Principles of Algebraic Geometry", John 
Wiley and Sons, New York, 1978

\bibitem[Gr90]{Gr90} Gromov, M., {\it Convex sets and K\"ahler manifolds}, Advances in 
Differential Geometry and Topology, F. Triceri ed., World Scientific, 
Teaneck, NJ, 1990; 1--38 

\bibitem[GS82]{GS82}  Guillemin V., and S. Sternberg, {\it Convexity properties of the
moment mapping I}, Invent. Math {\bf 67}(1982), 491--513

\bibitem[GS84]{GS84} ---, ``Symplectic 
Techniques in Physics,'' Cambridge Univ.\ Press, 1984

\bibitem[HS10]{HS10} Heinzner, P., and P.~Sch\"utzdeller, {\it 
Convexity properties of gradient maps},  Adv. Math. {\bf 225:3} (2010), 1119--1133

\bibitem[HNP94]{HNP94} Hilgert, J., K.-H. Neeb, and W. Plank, {\it Symplectic convexity 
theorems and coadjoint orbits}, Comp. Math.\ {\bf 94} (1994), 129--180 

\bibitem[Hi70]{Hi70} Hironaka, H., {\it Desingularization of complex-analytic varieties}, Actes Congr\`es intern. Math., Paris {\bf 2}(1970), 627--631

\bibitem[Ho84]{Ho84} Hovanskii, A., {\it Fewnomials and Pfaff manifolds}, I. C. M. 1983, 
Warszaw (1984), 549--565

\bibitem[Ju81]{Ju81} Jurkiewicz, J., {\it On the complex projective torus embeddings, the 
associated variety with corners and Morse functions}, Bull. Acad. Polon. Sci. 
Ser. Sci. Math. {\bf 29}(1981), 21--27

\bibitem[Kl13]{Kl13} Klartag, B., {\it Poincar\'e inequalities and moment maps}, 
Ann. Fac. Sci. Toulouse Math. (6) {\bf 22:1} (2013), 1--41 

\bibitem[LM87]{LM87} Libermann, P., and C. Marle, {\it Symplectic geometry 
and analytical mechanics}, Reidel, Dordrecht, 1987

\bibitem[Mi08]{Mi08} Mitter, S. K., {\it Convex optimization in infinite dimensional 
spaces}, in ``Recent Advances in Learning and Control, V.~D.~Blondel et al (eds.), 
LNCIS {\bf 371} (2008), 161--179

\bibitem[Ne92]{Ne92} Neeb, K. -- H., {\it Toric varieties and algebraic monoids}, Seminar 
Sophus Lie {\bf 2:2}(1992), 159--188

\bibitem[Ne98]{Ne98} ---, {\it On the complex and convex geometry of Ol'shanski\u\i{}
semigroups}, Annales de l'Inst.\ Fourier {\bf 48:1} (1998), 149--203 

\bibitem[Ne99]{Ne99} ---, {\it On the complex geometry of invariant domains 
in complexified symmetric spaces}, 
Annales de l'Inst.\ Fourier {\bf 49:1} (1999), 177--225

\bibitem[Ne00]{Ne00} ---, ``Holomorphy and Convexity in Lie Theory,''
Expositions in Mathematics {\bf 28}, de Gruyter Verlag, Berlin, 2000  

\bibitem[Od88]{Od88} Oda, T., ``Convex Bodies and Algebraic Geometry,'' 
Springer,  Berlin, Ergebnisse der Math. {\bf 15}, 1988

\bibitem[Od91]{Od91} ---, ``Geometry of Toric Varieties,'' Proc. of the Hyderabad Conf. on 
algebraic groups, S. Ramanan ed., Publ. for National Board for Higher Math. by 
Manoj Prakashar, 1991

\bibitem[Ok96]{Ok96} Okounkov, A., {\it Brunn--Minkowski inequalities for multiplicities}, Invent. Math. {\bf 125} (1996), 405--411

\bibitem[Ro70]{Ro70} Rockafellar, R. T., ``Convex Analysis,'' Princeton, New Jersey, 
Princeton University Press, 1970

\bibitem[Ro70]{Ro74} Rockafellar, R. T., {\it Conjugate duality and optimization}, 
in CBMS Series {\bf 16}, SIAM Publications 1974, 1--74 

\bibitem[Ru73]{Ru73} Rudin, W., ``Functional Analysis,'' McGraw Hill, 1973 

\bibitem[Te82]{Te82} Tessier, B., {\it Bonnesen-type inequalities in algebraic geometry}, in 
S. T. Yau ed., Princeton University Press, Annals of Math. Studies {\bf 
102}(1982), 85--105


\end{thebibliography}
\end{document}